\documentclass[twoside,onecolumn]{article}

\usepackage{blindtext} % Package to generate dummy text throughout this template

\usepackage[sc]{mathpazo} % Use the Palatino font
\usepackage[T1]{fontenc} % Use 8-bit encoding that has 256 glyphs
\usepackage{microtype} % Slightly tweak font spacing for aesthetics

\usepackage[english]{babel} % Language hyphenation and typographical rules

\usepackage[hmarginratio=1:1,top=32mm,columnsep=20pt]{geometry} % Document margins
\usepackage[hang, small,labelfont=bf,up,textfont=it,up]{caption} % Custom captions under/above floats in tables or figures
\usepackage{booktabs} % Horizontal rules in tables
\usepackage{graphicx}
\usepackage{lettrine} % The lettrine is the first enlarged letter at the beginning of the text

\usepackage{enumitem} % Customized lists
\setlist[itemize]{noitemsep} % Make itemize lists more compact

\usepackage{abstract} % Allows abstract customization
\usepackage{titlesec} % Allows customization of titles
\titleformat{\section}[block]{\large\scshape\centering}{\thesection.}{1em}{} % Change the look of the section titles
\titleformat{\subsection}[block]{\large}{\thesubsection.}{1em}{} % Change the look of the section titles

\usepackage{fancyhdr} % Headers and footers
\usepackage{titling} % Customizing the title section

\usepackage{hyperref} % For hyperlinks in the PDF
\usepackage{amsmath}
\usepackage{subfig}

\newcommand{\udisc}{\ensuremath{\mathbf{u}}}
\newcommand{\xdisc}{\ensuremath{\mathbf{x}}}

\usepackage{siunitx}
\usepackage{tabularx}
\usepackage{algorithm}
\usepackage{algpseudocode}
\usepackage{pifont}
\usepackage{multirow}
\newcommand{\xmark}{\ding{55}}%
\usepackage{booktabs}
\usepackage[most]{tcolorbox}
\pretitle{\begin{center}\Huge\bfseries} % Article title formatting
\posttitle{\end{center}} % Article title closing formatting
\title{A Nonlinear Physics-based Reduced Order Model with Convolutional-based Operator Compression} % Article title
\author{Anna Ivagnes  \\ \small SISSA, International School for Advanced Studies, \\ \small Mathematics Area, mathLab, Trieste, Italy. \\ \small  \href{mailto:aivagnes@sissa.it}{aivagnes@sissa.it} \normalsize \and Giovanni Stabile \\ \small Sant'Anna School of Advanced Studies
\\ \small  The Biorobotics Institute, %V.le R. Piaggio 34, 56025,
Pontedera, Pisa, Italy. \\ \small  \href{mailto:giovanni.stabile@santannapisa.it}{giovanni.stabile@santannapisa.it} \normalsize \and Gianluigi Rozza  \\ \small  SISSA, International School for Advanced Studies, \\\small  Mathematics Area, mathLab, Trieste, Italy. \\ \small \href{mailto:grozza@sissa.it}{grozza@sissa.it}}

\date{} % Leave empty to omit a date
\renewcommand{\maketitlehookd}{%
\begin{abstract}
\noindent Reduced-order models (ROMs) are widely used to accelerate the solution of parametrized partial differential equations, with projection-based methods relying on a linear representation of the governing equations in a low-dimensional space. Their efficiency, however, is fundamentally limited by the Kolmogorov $N$-width of the solution manifold: when the $N$-width decays slowly, accurate approximations require a large number of linear modes. This limitation has motivated nonlinear ROMs, which typically replace the linear representation of the solution field with nonlinear latent representations learned through autoencoder architectures. In this work, we propose a different perspective by introducing nonlinearity directly into the representation of the governing operators. Rather than compressing only the solution field, we learn nonlinear low-dimensional representations of the full-order operators while retaining an explicit reduced system of equations. Convolutional architectures are employed to exploit the locality and spatial structure of differential operators, and the resulting compressed operators are coupled with a nonlinear representation of the solution through a training strategy that accounts for both reconstruction accuracy and the error of the reduced-order solution. To address large sparse operators arising from fine and unstructured discretizations, continuous convolutions are used to operate directly on their non-zero entries. Finally, radial basis function interpolation enables the prediction of the compressed operators at unseen parameter configurations, providing an efficient online stage without assembling the corresponding full-order system. Numerical experiments on three parametrized problems demonstrate the potential of nonlinear operator compression as an alternative strategy for equation-based ROMs.
\end{abstract}
}

\begin{document}

% Print the title
\maketitle

%----------------------------------------------------------------------------------------
%	ARTICLE CONTENTS
%----------------------------------------------------------------------------------------

\section{Introduction}
\label{sec:intro}
Reduced-order models (ROMs) provide an effective framework for reducing the computational cost of high-fidelity simulations of parametrized partial differential equations (PDEs), while retaining the dominant features of the underlying physical system. Among the most established approaches, projection-based methods such as Proper Orthogonal Decomposition (POD)-Galerkin approaches construct a low-dimensional approximation of the solution by projecting the full-order problem onto a linear reduced space. Despite their mathematical robustness and computational efficiency, these approaches are fundamentally limited by the ability of a linear subspace to represent the solution manifold. In particular, their approximation properties are closely related to the Kolmogorov $n$-width of the solution manifold~\cite{cohen2015approximation, maday2002reduced, pinkus2012n}, which may decay slowly for nonlinear or transport-dominated problems. In such cases, an accurate representation requires a large number of basis functions, compromising the efficiency of the resulting ROM.
These limitations have motivated the development of \textbf{nonlinear model reduction techniques} capable of representing solution manifolds beyond linear subspaces. Among them, neural-network-based approaches, and in particular autoencoders (AEs), have received increasing attention~\cite{hinton2006reducing}. An autoencoder learns a nonlinear mapping between the high-dimensional solution space and a low-dimensional latent space through an encoder and a decoder, thereby providing a flexible nonlinear representation of the solution manifold. Several works have combined autoencoders with reduced-order modelling, using either parameter-to-latent mappings (in a \emph{nonintrusive} way)~\cite{lee2020model, pichi2024graph, han2022predicting, fresca2020deep, fresca2021comprehensive, franco2023deep} or learned latent dynamics to obtain an online approximation of the solution~\cite{vlachas2018data, otto2019linearly, kani2017dr, conti2023reduced, mucke2021reduced, rojas2021reduced, lehtimaki2022accelerating, linot2022data, farenga2025latent, farenga2022neural, tomada2026latent}. Such approaches have demonstrated the potential of nonlinear representations to achieve substantially lower-dimensional models than classical linear projection methods.
Despite these advances, most existing nonlinear ROMs based on neural networks focus primarily on the direct approximation or evolution of the latent variable. In this setting, the governing equations and their underlying algebraic structure are generally not explicitly retained in the learned reduced representation. For instance, once the latent variables have been obtained, their evolution may be entirely determined by a learned nonlinear map or neural differential equation, effectively replacing the original numerical solver with a data-driven surrogate. Although this paradigm can provide accurate predictions, it may reduce the interpretability of the reduced model and make it more difficult to preserve structural properties inherited from the full-order discretization.

This observation motivates a different perspective on nonlinear ROMs: \emph{rather than learning a nonlinear representation of the solution itself, can nonlinear mappings be used to compress the operators defining the governing equations?} In a classical projection-based ROM, the reduced operators are obtained by projecting the full-order operators onto a linear reduced space.
%For a full-order operator $\mathbf{A}\in\mathbb{R}^{N\times N}$ and a reduced basis $\mathbf{V}\in\mathbb{R}^{N\times r}$, for example, the corresponding reduced operator is given by
%\[ \mathbf{A}_r=\mathbf{V}^T\mathbf{A}\mathbf{V}.\]
The reduced problem therefore retains the algebraic structure of the original discretized equations, while operating in a space of reduced dimension. We propose to generalize this concept by replacing the linear projection with a nonlinear, trainable mapping that directly compresses the full-order operators into low-dimensional representations.
The proposed approach can be interpreted as a \textbf{nonlinear counterpart of projection-based model reduction}, with the following key novelties.
\begin{itemize}
    \item First, we introduce a \textbf{nonlinear operator compression strategy}, in which nonlinear neural architectures are employed to compress the full-order differential operators rather than the solution field itself. In this work, convolutional architectures are considered to exploit the spatial locality of the operators while providing a nonlinear generalization of the linear projection employed in classical projection-based ROMs.
    \item Second, we investigate the use of \textbf{continuous convolutional architectures}~\cite{coscia2023continuous} for operator compression. Differential operators arising from PDE discretizations are typically sparse and exhibit local connectivity, properties that can be naturally exploited by continuous convolutions. In contrast to standard convolutions, continuous convolutions can operate directly on unstructured spatial discretizations, avoiding the need to map the operators onto a structured grid and allowing their sparse representation to be processed directly.
    \item Finally, the proposed framework leads to a \textbf{nonlinear equation-based ROM}. The compressed operators are used to formulate and solve a reduced system of equations, while the nonlinear decoder reconstructs the corresponding full-order solution. In this way, the proposed methodology combines the nonlinear representation capabilities of neural networks with the equation-based nature of classical reduced-order modelling: the neural network does not directly replace the governing equations, but instead provides a nonlinear mechanism to construct their reduced representation.
\end{itemize}
The remainder of the manuscript is organized as follows. Section \ref{sec:methods} briefly recalls the background methodology on convolutional and continuous convolutional architectures, while Section \ref{sec:ceae} presents the nonlinear ROM introduced in this work. Section \ref{sec:results} investigates the performance of the method in three two-dimensional parametrized test cases: a Poisson case on a square, an advection-diffusion test case on a square, and a time-dependent viscous Burgers case in a backward facing step domain. A discussion on the computational time and resources needed for the three test cases can be found in Subsection \ref{subsec:cpu-time}. Finally, Section \ref{sec:conclusions} presents the main conclusions of this work, with several possible research directions.

\section{Background methodology: convolutional architectures}
\label{sec:methods}
This section briefly recalls how convolutional and continuous convolutional layers operate.

\textbf{Convolutional neural networks (CNNs)} are a class of neural architectures specifically designed to process \emph{structured} data like images. In contrast to fully connected networks, which typically require the input to be represented as a one-dimensional vector, convolutional architectures preserve the spatial organization of the input and exploit local correlations. This is achieved through two main properties of convolutional layers: \emph{local connectivity}, whereby each output value depends only on a local region of the input, and \emph{weight sharing}, whereby the same set of learnable parameters is applied throughout the computational domain. These properties allow convolutional networks to extract spatial features while keeping the number of trainable parameters relatively independent of the overall size of the input.

Consider an input tensor
\begin{equation}
\mathbf{S}
\in
\mathbb{R}^{C_{\mathrm{in}}\times H_{\mathrm{in}}\times W_{\mathrm{in}}},
\end{equation}
where $C_{\mathrm{in}}$ denotes the number of input channels and $H_{\mathrm{in}}$ and $W_{\mathrm{in}}$ are the spatial dimensions. A convolutional layer is characterized by a set of learnable kernels
\begin{equation}
\mathbf{K}
\in
\mathbb{R}^{C_{\mathrm{out}}\times C_{\mathrm{in}}\times H_K\times W_K},
\end{equation}
where $C_{\mathrm{out}}$ is the number of output channels and $H_K\times W_K$ is the kernel size. Denoting by $\mathbf{S}^{(c)}$ the $c$-th input channel and by $\mathbf{K}^{(d,c)}$ the kernel connecting input channel $c$ to output channel $d$, the output of the convolutional layer can be written as
\begin{equation}
\mathbf{Y}^{(d)}
=
h\left(
\sum_{c=1}^{C_{\mathrm{in}}}
\left(
\mathbf{K}^{(d,c)} * \mathbf{S}^{(c)}
+ b^{(d)}
\right)
\right),
\label{eq:conv-layer}
\end{equation}
where $h$ is a nonlinear activation function and $b^{(d)}$ is the bias associated with the $d$-th output channel. In the discrete setting, the convolution operation is evaluated locally as
\begin{equation}
\left( \mathbf{K}^{(d, c)} \ast \mathbf{S}^{(c)} \right)_{i, j}=\sum_{p=1}^{H_K} \sum_{q=1}^{W_K} K_{p, q}^{(d, c)}\, S_{i+p, j+q}^{(c)}.
    \label{eq:convolution}
\end{equation}
The same kernel is then translated across the entire spatial domain, which is the source of the weight-sharing property.

\begin{figure}[htpb!]
    \centering
    \subfloat[Discrete convolution]{\includegraphics[width=0.39\linewidth, trim={14cm 10cm 14cm 10cm}, clip]{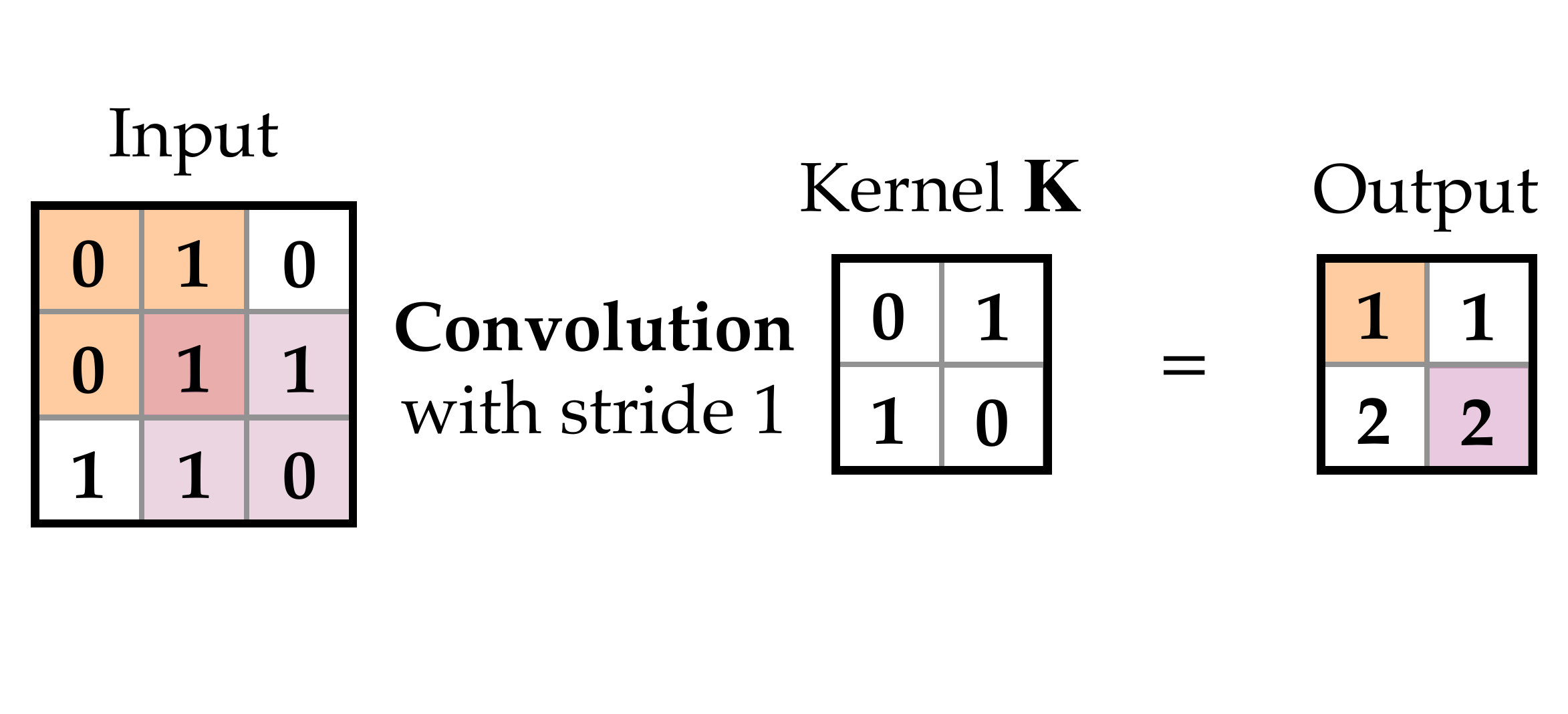} \label{subfig:discrete-conv}}
    \hfill
    \subfloat[Transposed discrete convolution]{\includegraphics[width=0.59\linewidth, trim={5cm 9cm 5cm 10cm}, clip]{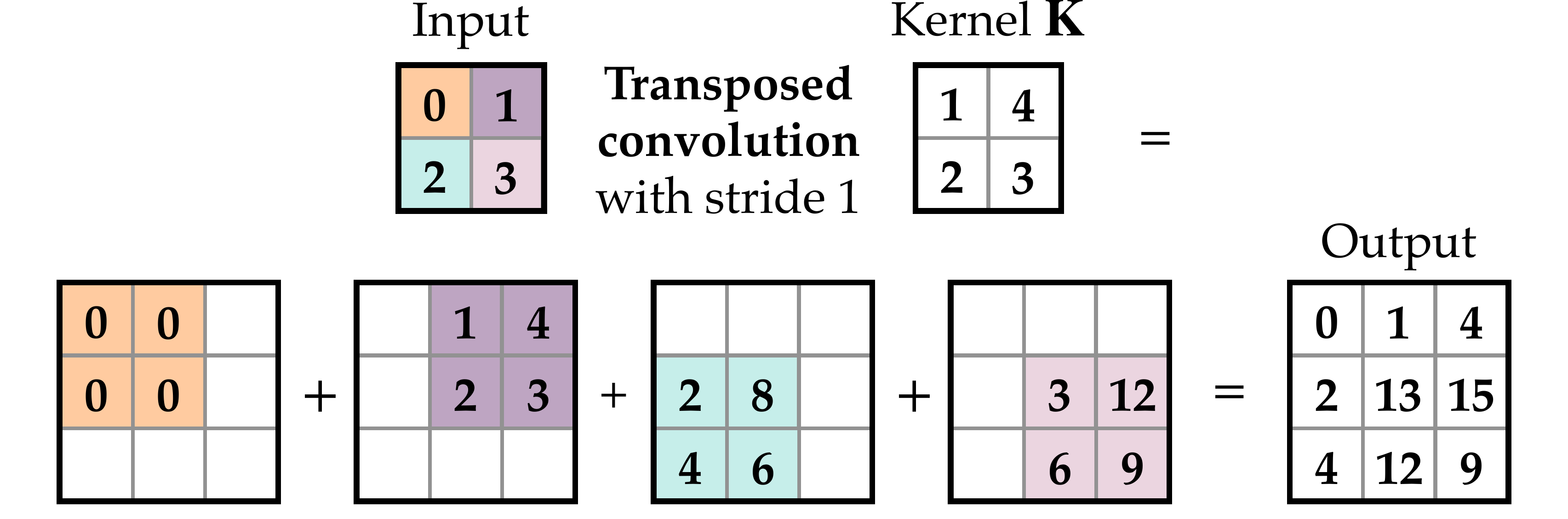} \label{subfig:discrete-transp-conv}}
    \caption{Examples of direct and transposed discrete convolution operations.}
    \label{fig:discrete-conv-operation}
\end{figure}

The spatial resolution of the feature maps can be controlled through additional convolutional hyperparameters, such as the stride, padding, and dilation. In particular, the stride determines the displacement of the kernel between consecutive applications and can therefore be used to reduce the spatial resolution. Padding controls the treatment of the boundary of the computational domain, while dilation increases the receptive field of the kernel without increasing its number of trainable parameters. An example of a convolutional and a transposed convolution operation is represented in Figure \ref{fig:discrete-conv-operation}.

By stacking convolutional layers and nonlinear activations, CNNs can consequently construct hierarchical representations in which progressively more complex spatial features are extracted.

Convolutional layers can also be incorporated into an autoencoder architecture. In a convolutional autoencoder (CAE), the encoder progressively maps the high-dimensional input to a low-dimensional latent representation, while the decoder reconstructs the original quantity from this representation. The decoder commonly employs transposed convolutional layers to perform learnable upsampling and recover the original spatial resolution. This provides a nonlinear dimensionality-reduction framework~\cite{lee2020model} that preserves the spatial structure of the input throughout the encoding and decoding operations.

\begin{figure}[htpb!]
    \centering
    \includegraphics[width=\linewidth, trim={0 13cm 0 9cm}, clip]{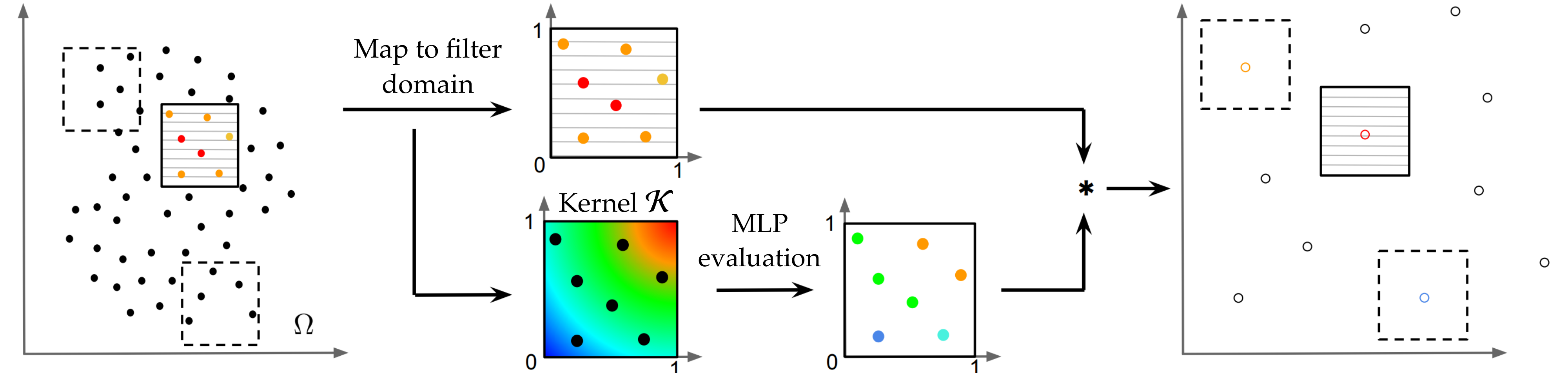}
    \caption{Sketch of the continuous convolutional operation~\cite{coscia2023continuous}: the unstructured points are mapped into the filter domain $\mathbb{X}\times \mathbb{Y}=[0, 1]^2$, then a MLP kernel $\mathcal{K}$ is used to approximate the filter values, and finally a convolution between the mapped input values and the filter values is computed.}
    \label{fig:continuous-conv}
\end{figure}

Despite these advantages, standard discrete convolutions are naturally formulated for data defined on structured grids. This requirement can be restrictive for computational problems involving unstructured spatial discretizations. A \textbf{continuous convolutional layer}~\cite{coscia2023continuous} addresses this limitation by replacing the discrete convolutional kernel with a learnable continuous function. For a continuous input field $\mathcal{S}$ and a continuous kernel $\mathcal{K}$, the convolution in a two-dimensional domain can be expressed as
\begin{equation}
\mathcal{Y}(x, y)=(\mathcal{S} \ast \mathcal{K})(\boldsymbol{x})=\int_{\mathbb{X}} \int_{\mathbb{Y}} \mathcal{S}(x + \tau_x, y+\tau_y) \cdot \mathcal{K}(\tau_x, \tau_y)\,d \tau_x\, d \tau_y,
    \label{eq:cont-convolution}
\end{equation}
where $\mathcal{X}\times\mathcal{Y}$ defines the local domain of the convolutional kernel. In contrast to the discrete case, the kernel is a continuous function of the relative coordinates. In practice, it can be represented by a trainable multilayer perceptron, while the integral in Equation~\eqref{eq:cont-convolution} is approximated through a discrete summation over the available input points\begin{equation}
\mathcal{Y}(\tilde{x}_i, \tilde{y}_i)=\sum_{x_i \in \mathbb{X}}\sum_{y_i \in \mathbb{Y}} \mathcal{S}(x_i + \tau_x, y_i+\tau+y) \cdot \mathcal{K}(x_i, y_i),
    \label{eq:cont-conv-sum}
\end{equation}
where $(\tau_x, \tau_y)$ are the strides positions, while $(\tilde{x}_i, \tilde{y}_i)$ are the points obtained by taking the centroid of the filter position mapped onto the $\Omega$ domain. A schematic illustration of the continuous convolutional layer is represented in Figure \ref{fig:continuous-conv}~\cite{coscia2023continuous}.

This formulation allows the convolutional operation to be applied directly to irregularly distributed data, without requiring the input field to be interpolated onto a regular grid~\cite{coscia2023continuous}. The continuous formulation can therefore be particularly advantageous for data arising from unstructured computational meshes.

The locality and weight-sharing properties of convolutional architectures also make them suitable for learning nonlinear representations of differential operators, represented by matrices, which are structured data. This observation provides the basis for the nonlinear ROM introduced in the following section.

\section{Nonlinear ROM based on operator compression}
\label{sec:ceae}
\begin{figure}
    \centering
    \includegraphics[width=\linewidth]{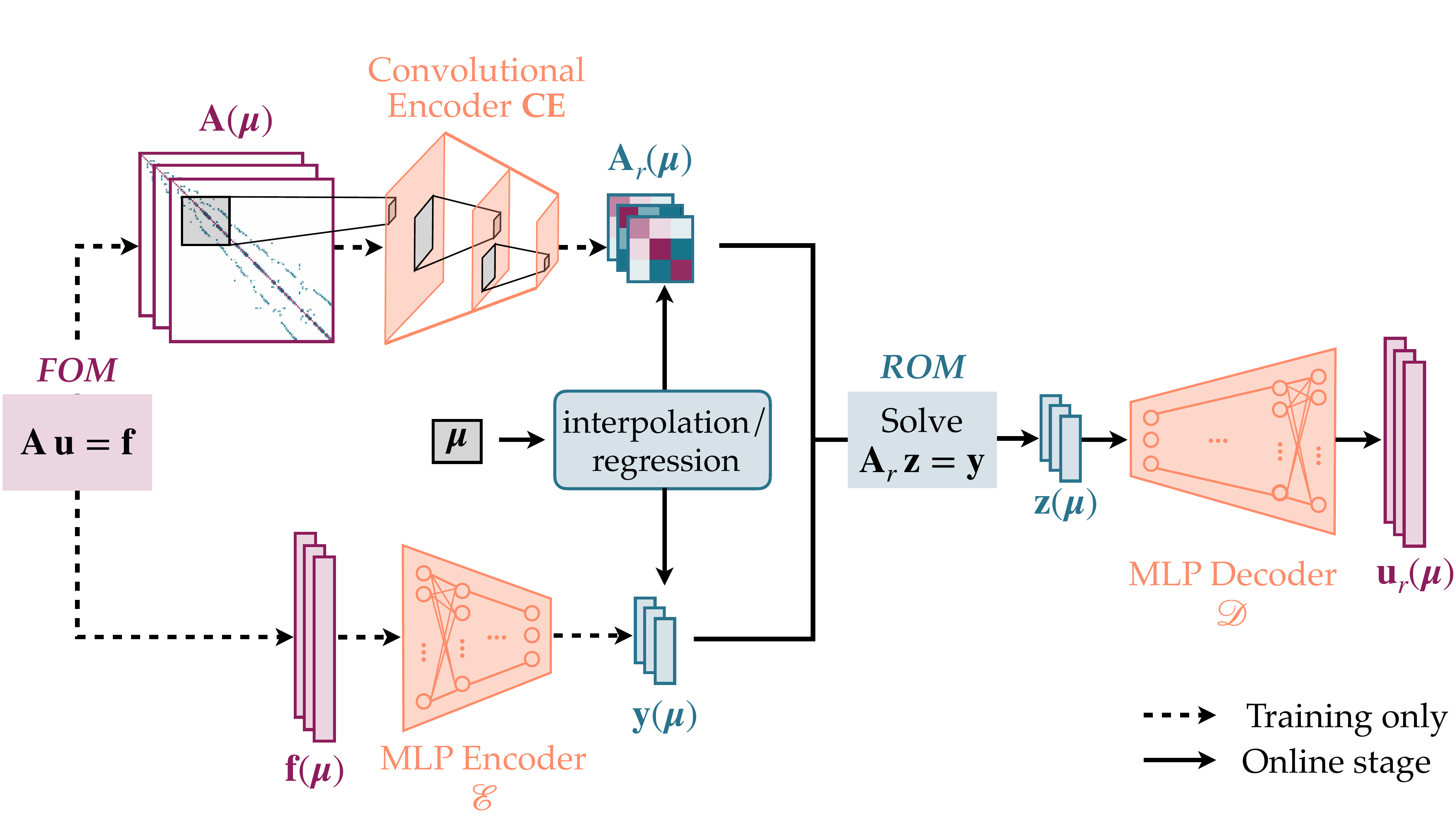}
    \caption{Sketch of the CE-AE pipeline: the convolutional encoder CE and the solution autoencoder AE are jointly trained offline to minimize the ROM residual, while an interpolation technique is adopted online to predict the parametrized reduced vectors and operators.}
    \label{fig:ceae-scheme}
\end{figure}

We present here the novel nonlinear ROM proposed in this work, which is sketched in Figure \ref{fig:ceae-scheme}.

Our goal is to propose a nonlinear version of the well-known projection-based intrusive ROM approaches, based on neural networks.

We start by considering a full-order system in the form $\mathbf{A} \mathbf{u}=\mathbf{f}$, where $\mathbf{A}\in \mathbb{R}^{N \times N}$, $\mathbf{u} \in \mathbb{R}^{N}$, and $\mathbf{f} \in \mathbb{R}^{N}$. We aim to compress the FOM into a reduced order system of prescribed dimension $r \ll N$. We seek three different mappings:
\begin{itemize}
    \item a nonlinear convolutional encoder \textbf{CE}: $\mathbf{A} \mapsto \mathbf{A_r}$, with $\mathbf{A_r}\in\mathbb{R}^{r \times r}$. To learn this nonlinear mapping, we employ CNNs. Originally developed for structured data such as images, CNNs are naturally suited to discrete differential operators, which can also be interpreted as structured inputs;
    \item a nonlinear encoder $\mathcal{E}: \mathbf{f} \mapsto \mathbf{y}$, with $\mathbf{y} \in \mathbb{R}^{r}$ is the latent representation of the full-order field $\mathbf{f}$;
    \item a nonlinear decoder $\mathcal{D}:\mathbf{z} \mapsto \mathbf{u}_r$, where $\mathbf{z}$ is found by solving the reduced $r$-dimensional system $\mathbf{A_r} \mathbf{z}=\mathbf{y}$.
    In a classical ROM fashion, the goal is to obtain a decoded solution $\mathbf{u}_r \in \mathbb{R}^{N}$ as close as possible to the full-order reference $\mathbf{u}$.
\end{itemize}

Both the matrix $\mathbf{A}$ and the vector-valued fields $\mathbf{u}$ and
$\mathbf{f}$ are in general parameter-dependent, namely $\mathbf{A}(\boldsymbol{\mu})$, $\mathbf{f}(\boldsymbol{\mu})$, and $\mathbf{u}(\boldsymbol{\mu})$.

We consider as encoder $\mathcal{E}$ and decoder $\mathcal{D}$ classical MLP neural networks, but other choices are possible, such as convolutional~\cite{lee2020model} or continuous convolutional neural networks~\cite{coscia2023continuous}.

Moreover, we consider the MLP encoder and decoder as parts of the same AE, e.g., having a symmetric architecture in the layers and activation functions, but different choices are possible. More in detail, about the AE architecture, we remark that the encoder and decoder weights are not symmetric but jointly trained without any structural constraint.

\begin{tcolorbox}[breakable,
    colback=gray!10,
    colframe=gray!10,
    boxrule=0pt,
    arc=2pt
]
\textbf{Remark: about the solution encoder.}\\
It is important to highlight that we compress the right-hand side term $\mathbf{f}$ with the same encoder used to compress the solution $\mathbf{u}$, as the two vectors formally belong to the same vector space. This is consistent with the linear ROM setting, where both the full-order solution and the right-hand side are represented in the same underlying discretized space. However, this assumption may introduce an excessive inductive bias in the training process. Indeed, although $\mathbf{u}$ and $\mathbf{f}$ belong to the same vector space, their spatial structures and relevant features can differ significantly. Enforcing a shared encoder therefore restricts the admissible latent space and may prevent the network from learning the most suitable reduced representation for the solution, negatively affecting the convergence.
\end{tcolorbox}

%loss function considered
In standard AEs, the loss function that is typically adopted is the discrepancy between the reconstructed solution, namely, adopting the notation of Figure \ref{fig:ceae-scheme}:
\begin{equation}
    \mathcal{L}_{\text{AE}}=\dfrac{1}{N_{\mu}}\sum_{i=1}^{N_{\mu}}\|\udisc(\boldsymbol{\mu}_i)-\udisc_r(\boldsymbol{\mu}_i)\|_2^2,
    \label{eq:loss-ae-simple}
\end{equation}
where $N_{\mu}$ is the number of training snapshots, and $\boldsymbol{\mu}$ is the solution set of parameters.

The novel architecture has instead a loss function in the form
\begin{equation}
    \mathcal{L}_{\text{CE-AE}} = \mathcal{L}_{\text{CE}}+\mathcal{L}_{\text{AE}},
    \label{eq:loss-ceae}
\end{equation}
where we distinguish the following contributions.

\begin{itemize}
    \item The \emph{intrusive} contribution $\mathcal{L}_{\text{CE}}$, namely the error made by solving the reduced system. In particular, it is computed as follows.
    \begin{enumerate}
        \item Compute the reduced operator $\mathbf{A}_r=\text{CE}(\mathbf{A})$, and the reduced vector $\mathbf{y}=\mathcal{E}(\mathbf{f})$;
        \item Solve the reduced system $\mathbf{A}_r \mathbf{z}=\mathbf{y}$ and find $\mathbf{z}$.
        \item Backmap the reduced solution $\mathbf{z}$ to the full-order space using the decoder $\mathcal{D}$, and compute the \emph{intrusive} loss as:
        \begin{equation}
            \mathcal{L}_{\text{CE}}=\dfrac{1}{N_{\mu}}\sum_{i=1}^{N_{\mu}}\|\mathbf{u}-\mathcal{D}(\mathbf{z})\|_2^2.
            \label{eq:loss-intrusive-cae}
        \end{equation}
    \end{enumerate}
    \item The \emph{nonintrusive} contribution $\mathcal{L}_{\text{AE}}$, representing the AE$=\mathcal{D}\circ \mathcal{E}$ reconstruction error solely. This contribution is the one of Equation \eqref{eq:loss-ae-simple}, with $\mathbf{u}_r=(\mathcal{D}\circ \mathcal{E})(\mathbf{u})$.
\end{itemize}

We call the proposed approach \textbf{CE-AE} (Convolutional Encoder-AutoEncoder).

\begin{tcolorbox}[breakable,
    colback=gray!10,
    colframe=gray!10,
    boxrule=0pt,
    arc=2pt
]
\textbf{Remark: about the compression of more complex FOM systems.}\\
The proposed logic may be extended to more complex systems by applying the nonlinear compression to each full-order operator appearing in the governing equations. In the presence of two or more operators, we can  modify our approach in two different ways: {\textbf{(i)}} consider different CEs, which separately compress the given operators; {\textbf{(ii)}} consider a unique CE with more input channels, each processing one operator.

More complex nonlinear FOMs in the form $\mathcal{R}(\mathbf{u};\xdisc, \boldsymbol{\mu})=\mathbf{0}$ lead to a ROM of the form $\mathcal{R}_r(\mathbf{z};\boldsymbol{\mu})=\mathbf{0}$. In this case, the exact resolution of the ROM system can be handled in two ways:
\begin{itemize}
    \item Iteratively solve the nonlinear system using a Newton method. In this case, the reduced solution $\mathbf{z}$ would be the \emph{exact} one, and the loss contribution remains the same as in Equation \eqref{eq:loss-intrusive-cae}.
    \item Replace the loss contribution of Equation \eqref{eq:loss-intrusive-cae} with the residual norm, leading to:
    \begin{equation}
        \mathcal{L}_{\text{CE}}=\dfrac{1}{N_{\mu}}\sum_{i=1}^{N_{\mu}}\|\mathcal{R}_r(\mathbf{z};\boldsymbol{\mu}_i)\|_2^2.
        \label{eq:loss-intrusive-nonlinear}
    \end{equation}
    This transforms the model in a PINN~\cite{raissi2019physics} at the reduced-order level, as done in previous works like~\cite{kim2022fast, chen2021physics, dave2025physics, hijazi2023pod}.
\end{itemize}
\end{tcolorbox}

In the case of parametrized operators $\mathbf{A}(\boldsymbol{\mu})$ and/or right-hand side vectors $\mathbf{f}(\boldsymbol{\mu})$, the previously described architecture cannot be directly employed for the evaluation of unseen parameter configurations $\boldsymbol{\mu}_{\star}$. Indeed, the proposed method requires access to the FOM quantities, namely the operators and/or the full-order vectors, which are not available at new parameter values.

Therefore, an additional parametrization step is required to enable prediction at unseen configurations. As illustrated in Figure \ref{fig:ceae-scheme}, we introduce a parametric mapping $\boldsymbol{\mu} \mapsto \mathbf{A}_r$ and/or $\boldsymbol{\mu} \mapsto \mathbf{y}$, which directly predicts the corresponding latent representation of the operators and vectors from the parameter space.

In this way, the reduced-order quantities can be efficiently reconstructed without requiring the assembly or evaluation of the underlying FOM system at the new configuration $\boldsymbol{\mu}_{\star}$.
In this case, we adopt a RBF interpolation approach for each element of the parametrized vector/matrix, but different approaches (like GPR or NN) are possible. We will indicate the proposed ROM, coupled with the RBF prediction at unseen states as \textbf{CE-AE-RBF}.

Additionally, in convection-dominated cases the full-order convection matrix is typically $\mathbf{C}(\mathbf{u})$, hence a mapping $\mathbf{z} \mapsto \mathbf{C}_r$ is here proposed, as we will see in the last test case considered (Section \ref{subsec:advec-cae}). In time-dependent test cases, this also allows the generalization of the approach to different time step sizes, as the time instance $t$ is not an explicit input of the reduced map.

\begin{tcolorbox}[colback=gray!10, breakable,
    colframe=gray!10,
    boxrule=0pt,
    arc=2pt]
    {\textbf{Remark: about the stability of the reduced system $\mathbf{A}_r \mathbf{z}=\mathbf{y}$.}}\\
We remark that in all the proposed convolutional approaches and in all the test cases, we consider as output matrix the sum of the encoder output matrix, and of the \textbf{stabilization term} $ \epsilon \mathbf{I}_r$, where $\mathbf{I}_r$ is the identity matrix of dimension $r \times r$ and $\epsilon=\num{1e-4}$.
The addition of this regularization term enhances the numerical stability and improves the conditioning of the reduced system $\mathbf{A}_r \mathbf{z}=\mathbf{y}$. This system is solved both in the offline stage, during the training process, and in the online stage, when the solution is predicted for unseen configurations.

Another possible strategy to prevent poor conditioning is to augment the loss function with a regularization term based on the condition number of the reduced operator. Nevertheless, the repeated computation of the condition number during training would significantly increase the computational cost of each optimization step. For this reason, we do not adopt this approach in the present work.
\end{tcolorbox}

A natural question is how to choose the reduced dimension $r$. In this work, we focus on very low-dimensional latent spaces ($r=2,3,4$), where the nonlinear representation is expected to provide the largest advantage over linear reduction techniques. For larger values of $r$, the accuracy of the POD approximation progressively improves, while the autoencoder performance tends to saturate.

\subsection{Continuous convolution for sparse operators: CCE-AE}
\label{subsec:cceae}
In real-world scenarios, the FOM system usually involves a very large number of degrees of freedom $N$, potentially reaching the order of millions. As a result, the corresponding FOM matrices may have prohibitively large dimensions ($N \times N$), preventing their direct manipulation and storage. This limitation can make the adoption of standard convolutional architectures impractical, as they often require memory resources that scale unfavourably with the input size.

However, it is worth noting that FOM operators are typically sparse, meaning that the majority of their entries are equal to zero.

The question now is: {\emph{How can we exploit this property to overcome the memory limitations associated with the direct treatment of full matrices?}}\\
In this work we consider the continuous convolutional architecture, specifically designed to operate on irregularly sampled data and can naturally exploit the sparsity structure of the operators.

\begin{figure}[htpb!]
    \centering
    \includegraphics[width=\linewidth, trim={3cm 22cm 3cm 0cm}, clip]{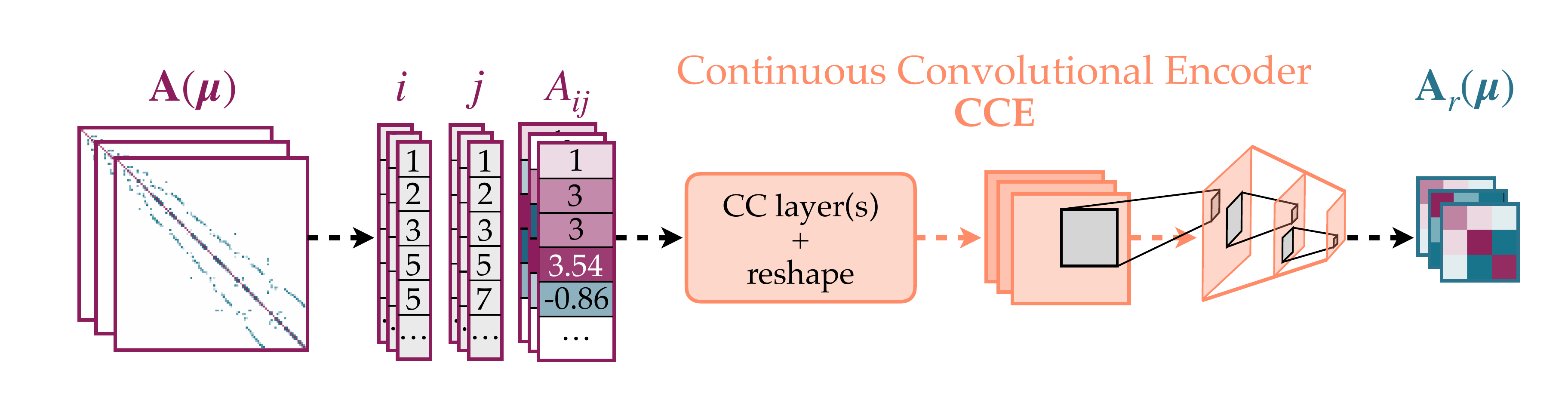}
    \caption{Sketch on the operator compression strategy using the continuous convolutional architecture, coupled with a standard convolutional network.}
    \label{fig:ccencoder}
\end{figure}

In the original paper on continuous convolutional layers~\cite{coscia2023continuous}, the input consists of the spatial coordinates (e.g., $(\mathrm{x},\mathrm{y})$ in a two-dimensional setting) together with the corresponding values of the physical field evaluated at those locations $\mathbf{u}(\mathrm{x},\mathrm{y})$. Following the same idea, we interpret the discrete operator as the realization of an underlying continuous field. Instead of providing the full $N \times N$ matrix as input, we consider only the entries corresponding to non-zero elements, providing their indices $(i,j)$ and the associated values. Consequently, the input dimensionality is reduced from $N \times N$ to $3 \times N_{\mathrm{nz}}$, where $N_{\mathrm{nz}} \ll N$ denotes the number of non-zero entries and is typically much smaller than $N$, depending on the structure of the considered operator. As represented in Figure \ref{fig:ccencoder}, one (or more) continuous convolutional layers are used to compress the original FOM operator to an intermediate dimension, which can be then further compressed by standard convolutional layers.

The nonlinear ROM approach using the continuous convolutional encoder architecture is here named CCE-AE. Analogously to the standard convolutional ROM (CE-AE), we refer to the efficient ROM version using the RBF interpolation as \textbf{CCE-AE-RBF}. This formulation offers two important advantages. First, it allows the proposed approach to be applied to different types of operators, such as mass, diffusion, or stiffness matrices, without requiring a modification of the architecture. Second, since the input is defined through the coordinates of the non-zero entries rather than through a fixed matrix representation, the method can naturally handle operators arising from different discretizations and computational grids.

Furthermore, if a CCAE is adopted also for the solution field, the entire framework could be extended to become \textbf{independent of the underlying discretization grid}. In this case, both the operator and the solution would be represented through continuous latent mappings, enabling the application of the reduced model across different grids and full-order discretizations. This extension is not investigated in the present work and is left as a direction for future research.

\section{Numerical results}
\label{sec:results}

Here we evaluate the performance of the proposed approach in three different test cases:
\begin{itemize}
    \item  a Poisson equation on a square with two parameters incorporated into the source term (Section \ref{subsec:poisson-cae});
    \item a parametrized advection-diffusion model defined on a square domain, with one parameter defining the diffusivity (Section \ref{subsec:advec-cae});
    \item a Burgers problem in a 2D backstep domain, where the time is the unique parameter considered (Section \ref{subsec:burgers-cae}).
\end{itemize}
For each test case, we assess the performance of the proposed CE-AE approach (and, for parametrized problems, of the CE-AE-RBF framework) by comparing it with standard linear reduction techniques, such as POD and POD-RBF, as well as nonlinear approaches based on AE and AE-RBF. We specify that we exploit the same AE architecture both in the CE-AE and in the AE solely, to guarantee a fair comparison.

The accuracy of the different methods is evaluated through the relative $L^2$ error of the reconstructed solution with respect to the corresponding FOM solution, together with a qualitative assessment of the obtained approximations. Furthermore, we provide a comparison of the loss behaviour of standard AEs and the proposed CE-AE approaches, investigating an ensemble of different NNs with the same architecture and different weights' initialization.

Finally, we analyse the reduced operators obtained through the proposed nonlinear compression strategy and discuss their underlying properties.

\subsection{Poisson test case}
\label{subsec:poisson-cae}
The first test case is a Poisson non-affine equation on the square domain $\Omega=[-1, 1]^2$, with two parameters $\boldsymbol{\mu}=(\mu_1, \mu_2) \in [-1, 1]^2$, written in the following form:
\begin{equation}
\begin{cases}
    \Delta u=g(\boldsymbol{x};\boldsymbol{\mu}) \text{ in }\Omega, \\
    u=0 \text{ on }\partial \Omega,
\end{cases}
\text{ with }g(\boldsymbol{x};\boldsymbol{\mu})=\exp{\left( -2 (x-\mu_1)^2 -2(y-\mu_2)^2\right)}.
    \label{eq:poisson-equation}
\end{equation}
We can then obtain the following weak formulation:
\[ a(u, v)=f(v;\boldsymbol{\mu}), \text{where}\,\,
a(u, v)=\int_{\Omega} \nabla u \cdot \nabla v \, d\boldsymbol{x}, \quad f(v;\boldsymbol{\mu})=\int_{\Omega} g(\boldsymbol{x},\boldsymbol{\mu}) \, d\boldsymbol{x},\]
for each test function $v \in H^1_0(\Omega)$. Using a FE discretization, the above formulation can be written as a linear system with parametrized right-hand side:
\begin{equation} \mathbf{A}\mathrm{u}(\boldsymbol{\mu}) = \mathrm{f}(\boldsymbol{\mu}),
\label{eq:poisson-fom-system}
\end{equation}
with $\mathbf{A} \in \mathbb{R}^{N\times N}$, $\mathrm{f} \in \mathbb{R}^{N}$, and $N=638$. The FE grid has been built using the \textit{gmsh} utility, and the FE simulation is performed using FEniCS~\cite{logg2010dolfin, logg2010dolfin, alnaes2015fenics}.

\begin{figure}[htpb!]
    \centering
    \subfloat[Parameters' space]{\includegraphics[width=0.4\linewidth]{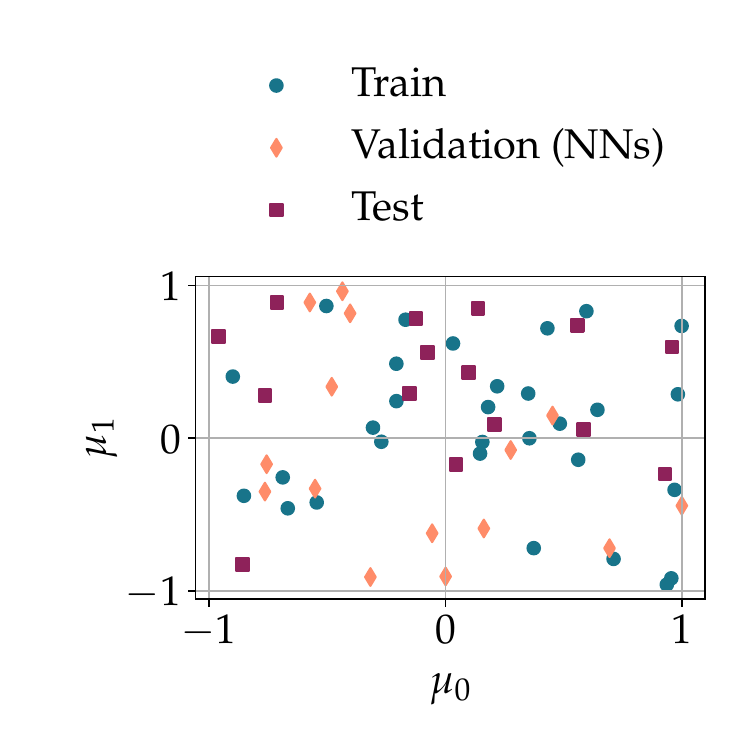}}
    \subfloat[POD singular values' decay]{\includegraphics[width=0.52\linewidth]{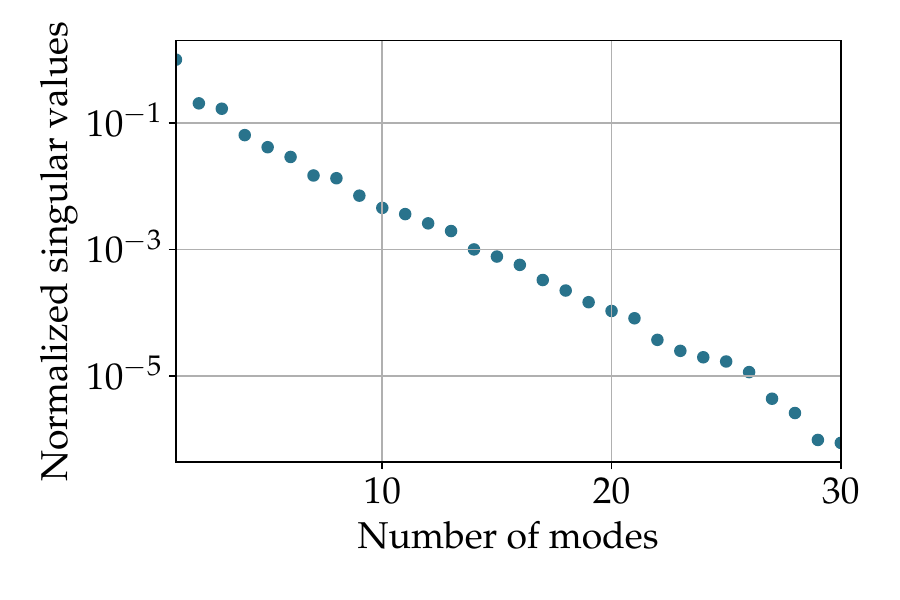}}
    \caption{\emph{Poisson test case.} Train and test parameters (left), and the POD singular values decay (right).}
    \label{fig:poisson-params}
\end{figure}

We collect a set of $30$ train and $15$ test snapshots, and then perform the POD considering the training data. Figure \ref{fig:poisson-params} represents the values of the train, test, and validation parameters, together with the POD singular values decay. The validation set, composed of $15$ snapshots, is used during the training of the NNs solely to monitor the different components of the loss function and assess the model's generalization performance.

Although the singular values exhibit an exponential decay, the decay rate is relatively slow, indicating that a large number of POD modes is required to capture a significant fraction of the solution energy.

\begin{figure}[htpb!]
    \centering
    \includegraphics[width=0.4\linewidth, trim={2cm 0 1cm 0}, clip]{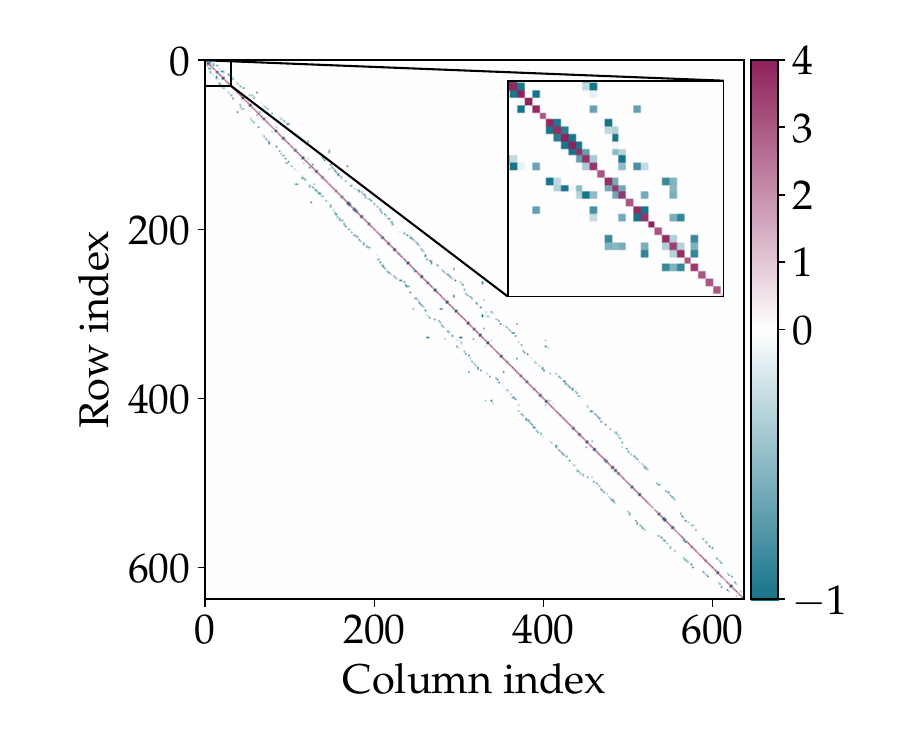}
    \caption{\emph{Poisson test case.} Full-order matrix $\mathbf{A}$.}
    \label{fig:poisson-full-matrix}
\end{figure}

Our objective is to achieve accurate ROM predictions in \emph{under-resolved} or \emph{marginally-resolved} regimes, namely by considering a reduced number of degrees of freedom.
For this reason, we analyse the results obtained considering $r=2, 3$, and $4$.
The diffusion matrix $\mathbf{A}$ in Equation \eqref{eq:poisson-fom-system}, whose sparsity pattern is represented in Figure \ref{fig:poisson-full-matrix}, is symmetric and positive-definite, namely $\mathbf{A}^T=\mathbf{A}$ and $\mathrm{x}^T \mathbf{A} \mathrm{x} > 0$, $\forall \mathrm{x} \neq \mathbf{0}$.

In the classical reduced basis framework, the corresponding reduced matrix $\mathbf{A}_r=\mathbf{V}^T \mathbf{A} \mathbf{V}$ preserves the symmetry and positive-definiteness of the FOM matrix $\mathbf{A}$, because the discrete modes matrix $\mathbf{V}$ is an orthogonal matrix. We then pose ourselves the following question:\\{\emph{Is the \textbf{optimal} reduced matrix $\mathbf{A}_r$ symmetric positive-definite?}}\\
To numerically answer this question, we investigate two compression variants, where the final reduced matrix is expressed in two different ways:
\begin{itemize}
    \item $\mathbf{A}_r = \mathbf{A}^{\prime}_r + \epsilon \mathbf{I}_r$, where $\mathbf{A}^{\prime}_r$ is the encoder output, and $\epsilon = \num{1e-4}$. The stabilization term $\epsilon \mathbf{I}_r$ improves numerical stability and conditioning when solving the reduced system. The resulting model is simply referred to as CE;
    \item $\mathbf{A}_r = (\mathbf{A}^{\prime}_r)^T \mathbf{A}^{\prime}_r+\epsilon \mathbf{I}_r$, which is by construction symmetric and positive-definite. This model is called S-CE.
\end{itemize}

\begin{figure}[htpb!]
    \centering
    \subfloat[Latent dimension $r=2$]{\includegraphics[width=\linewidth]{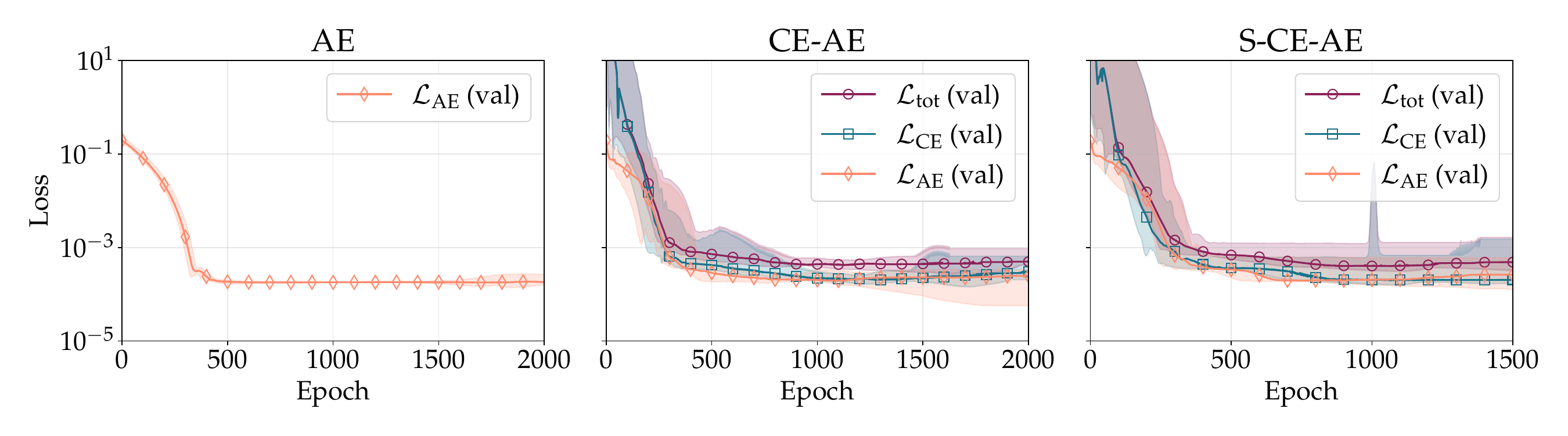}}\\
    \subfloat[Latent dimension $r=3$]{\includegraphics[width=\linewidth]{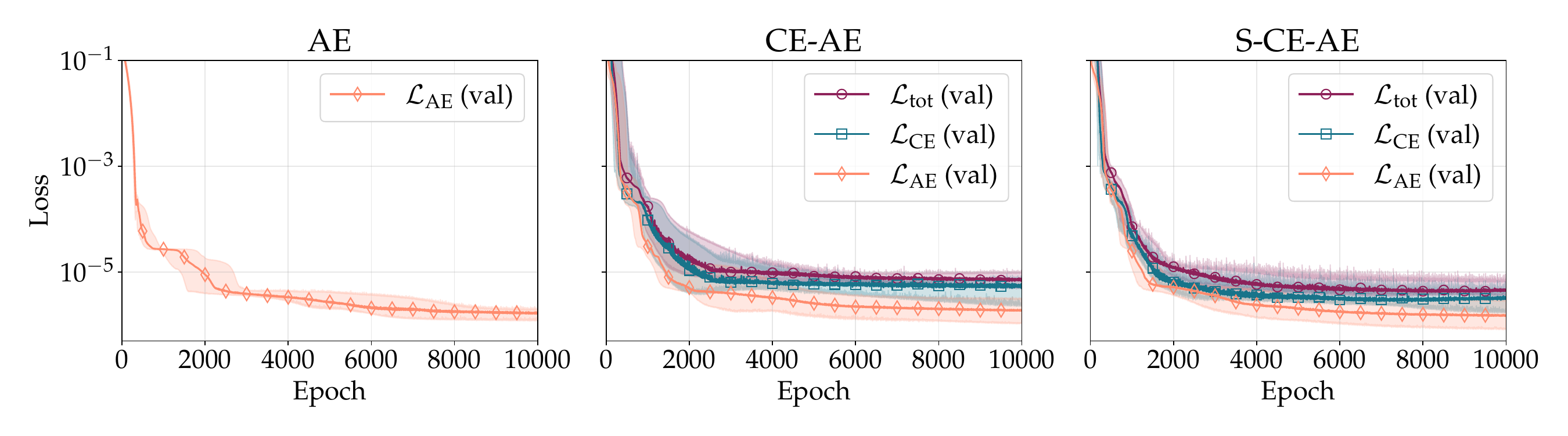}}\\
    \subfloat[Latent dimension $r=4$]{\includegraphics[width=\linewidth]{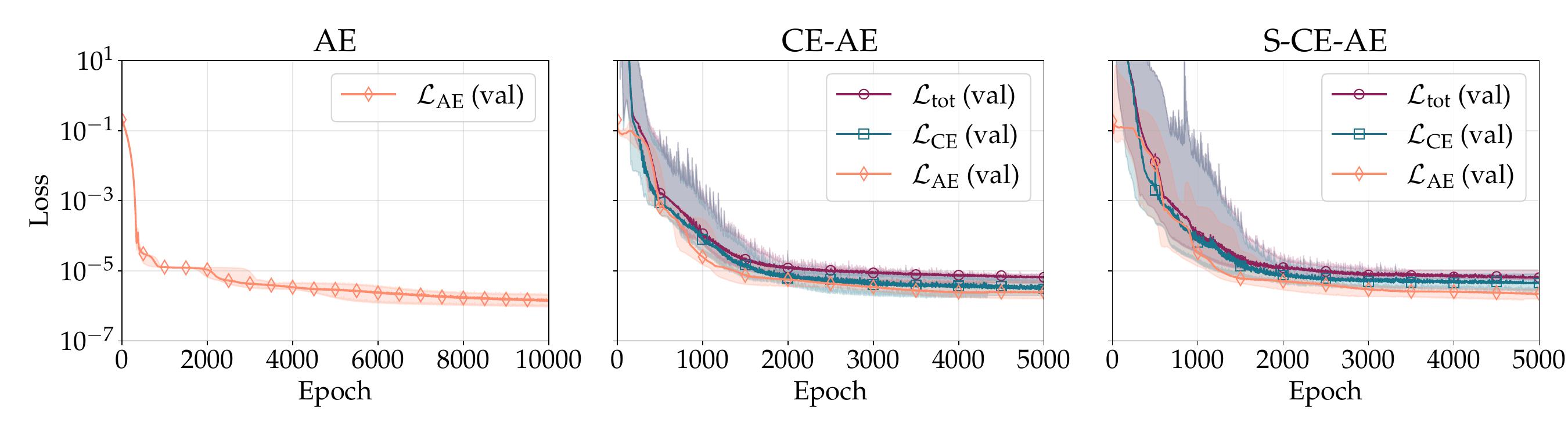}}
    \caption{\emph{Poisson test case.} Loss trend across epochs for a standard AE, and for the novel models CE-AE and S-CE-AE, at different latent dimensions, namely $r=2,3,$ and $4$.}
    \label{fig:poisson-losses}
\end{figure}

Before discussing the results, we first analyse the training behaviour of the different models by comparing the validation loss evolution of the AE, CE-AE, and S-CE-AE approaches. In particular, Figure \ref{fig:poisson-losses} reports the loss trends over the training epochs for $10$ neural networks sharing the same architecture but initialized with different random weights. The curves are presented in terms of their average value and associated confidence interval, in order to assess both the convergence behaviour and the variability induced by the initialization.

For all the considered latent dimensions $r=2, 3, 4$, we observe that the standard AE training exhibits, in general, a lower variability compared to the proposed CE-AE and S-CE-AE approaches. This behaviour is expected, since the latter models introduce an additional loss contribution, thus imposing a stronger constraint on the optimization process and potentially increasing the sensitivity to the initial conditions.

Moreover, the average value of the loss component associated with the AE reconstruction, $\mathcal{L}_{\text{AE}}$, at convergence is comparable among the three approaches for all the considered latent dimensions. However, for $r=2$ and $r=3$, the minimum values achieved during training, represented by the lower bound of the confidence intervals, are lower for both the CE-AE and S-CE-AE models. This suggests that the inclusion of the additional loss term $\mathcal{L}_{\text{CE}}$ can improve the convergence of the autoencoder reconstruction loss, allowing the network to reach lower reconstruction errors compared to the standard AE training alone.

\begin{figure}[htpb!]
    \centering
    \includegraphics[width=\linewidth, trim={3cm 0 7cm 0}, clip]{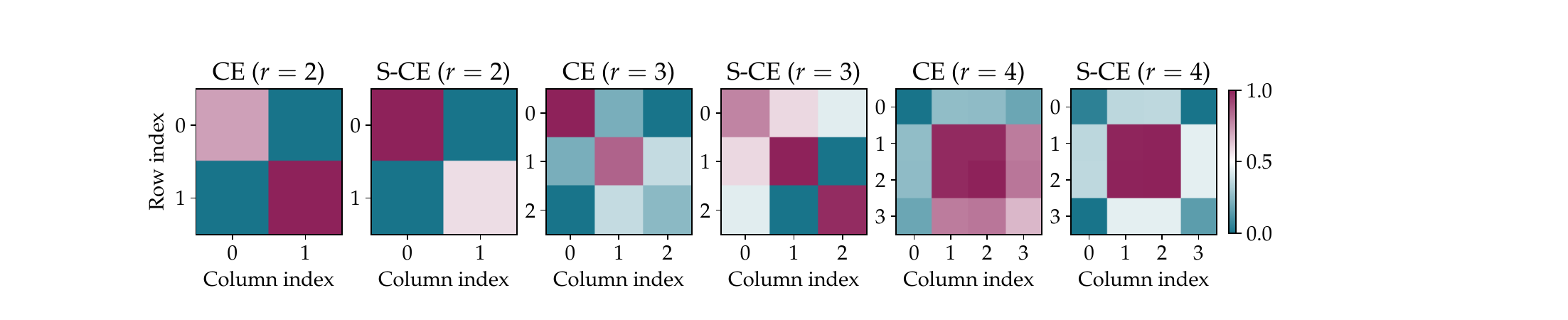}
    \caption{\emph{Poisson test case.} Reduced matrices $\mathbf{A}_r$, normalized in range $[0, 1]$. The matrices are obtained by compressing the full-order discrete matrix $\mathbf{A}$ with the convolutional encoder CE and with its symmetric version S-CE, at different latent dimensions $r=2, 3,$ and $4$.}
    \label{fig:poisson-matrices}
\end{figure}

The reduced matrices obtained with the proposed approaches are graphically represented in Figure \ref{fig:poisson-matrices}. The matrices are normalized between 0 and 1 to enable a qualitative comparison of their structures.

We first consider the case with latent dimension $r=2$. The reduced matrices obtained by the CE and S-CE models at convergence exhibit a similar structure, up to a permutation of the reduced degrees of freedom, which in this case corresponds to a transposition of the matrix. This behaviour can be attributed to the lack of an intrinsic ordering of the latent variables learned by the AE: different training runs may identify equivalent reduced representations associated with different latent coordinate arrangements.

Interestingly, both reduced operators are symmetric, despite the fact that no symmetry constraint is explicitly enforced in the CE architecture. The same behaviour is observed also for larger latent dimensions ($r=3$ and $r=4$), suggesting that the optimization process naturally drives the network towards symmetric reduced operators. This indicates that the optimal reduced representation identified by the neural network tends to preserve the structural properties of the original FOM operator, similarly to standard projection-based linear reduction techniques.

A similar correspondence between the global structures of the reduced operators can also be observed for the case $r=4$, further confirming the consistency of the learned representations across different latent dimensions.

\begin{figure}[htpb!]
    \centering
    \includegraphics[width=\linewidth]{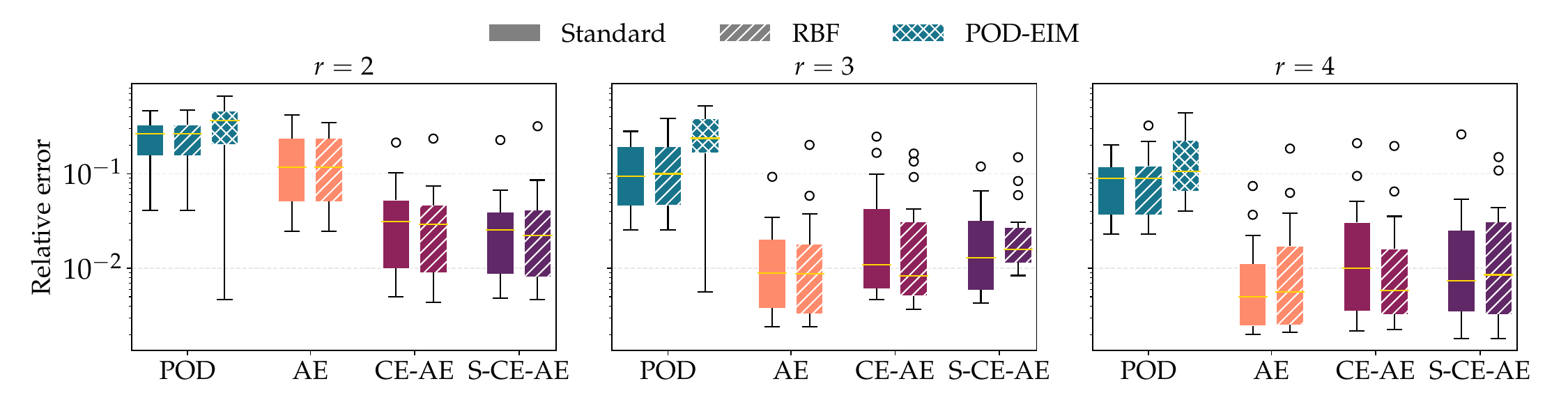}
    \caption{\emph{Poisson test case.} Relative $L^2$ test errors of different approximations with respect to the FOM reference, at the latent dimensions $r=2, 3$, and $4$.}
    \label{fig:poisson-boxplot}
\end{figure}

The relative $L^2$ test errors obtained by the proposed approaches, together with the corresponding POD-based results, are reported in Figure \ref{fig:poisson-boxplot}. The figure shows the box plots of the relative $L^2$ errors obtained using the best-performing network among the different training realizations considered in Figure \ref{fig:poisson-losses}. The selected network is then evaluated over a set of unseen parameter configurations (represented in Figure \ref{fig:poisson-params}).

For the POD framework, we consider three different error contributions: the reconstruction error (obtained by directly compressing and reconstructing the solution field), the POD-RBF error, and the POD-EIM error.

The POD-EIM approach is a projection-based intrusive method in which the reduced quantities (the right-hand side vector $\mathrm{y}$ in the present case) at unseen parameter configurations are approximated through an affine expansion of pre-computed EIM interpolation modes. In this work, we consider $r_{\text{EIM}}=2r$ interpolation modes, although different choices are possible depending on the desired accuracy and computational cost.

\begin{figure}[htpb!]
    \centering
    \includegraphics[width=\linewidth, trim={0 2cm 0 1cm}, clip]{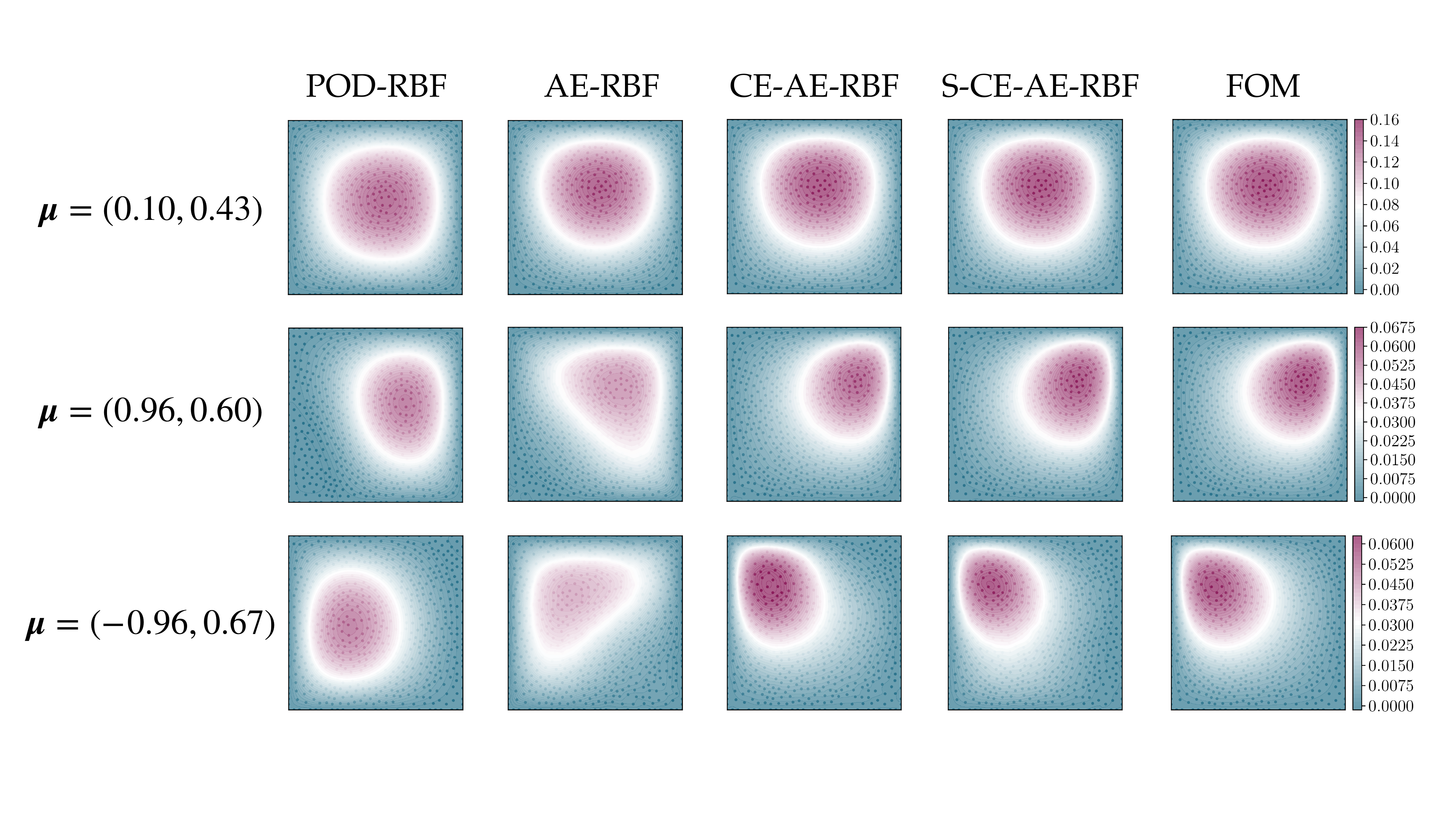}
    \caption{\emph{Poisson test case.} Approximated solutions for POD-RBF, AE-RBF, CE-AE-RBF, and S-CE-AE-RBF, compared with the FOM reference, at latent dimension $r=2$, and considering three test parameters.}
    \label{fig:poisson-fields}
\end{figure}

From the plots, we can draw the following conclusions.
\begin{itemize}
    \item For $r=2$, the proposed approaches achieve an accuracy improvement of approximately one order of magnitude compared to the standard AE. This result is further supported by the qualitative comparison reported in Figure \ref{fig:poisson-fields}, which shows that both the POD and the standard AE fail to accurately generalize the solution prediction to unseen configurations.
    \item At $r=3$ and $r=4$, the proposed approaches achieve results comparable to the AE-RBF.
\end{itemize}

As already discussed when analysing the loss behaviour, the inclusion of the convolutional encoder and the additional loss contribution does not negatively affect the performance of the autoencoder. Instead, the proposed formulation provides an additional regularizing effect, which stabilizes the optimization process and allows the network to reach comparable or even improved accuracy with respect to the standard AE.

\subsection{Advection-diffusion test case}
\label{subsec:advec-cae}

Here, we are interested in a parametrized advection-diffusion problem, defined on the unit square $\Omega=[0, 1]^2$, written as follows:

\begin{equation}
\begin{cases}
    -10^{-\mu} \Delta u + \beta \cdot \nabla u=1 \text{ in }\Omega,\\
    u=0 \text{ on }\partial \Omega,
\end{cases}
    \label{eq:adv-diff-model}
\end{equation}
where $\beta=(1, 1)$ is the advection velocity field, while the parameter $\mu$ determines the Péclet number of the problem. Increasing $\mu$ reduces the diffusion coefficient, leading to advection-dominated regimes, whereas smaller values of $\mu$ correspond to diffusion-dominated ones. The corresponding weak form can be written as:
\[
a(u, v;\mu)=f(v), \, \forall v \in H^1_0(\Omega),
\]
with \[
a(u, v;\mu)=\int_{\Omega}  10^{-\mu} \nabla u \cdot \nabla v\, d\boldsymbol{x} + \int_{\Omega}(\beta \cdot \nabla u)v \, d\boldsymbol{x}, \quad f(v)=\int_{\Omega}v \,d\boldsymbol{x}.
\]
Since the above formulation might become unstable at some $\mu$ values, we apply the SUPG stabilization, obtaining the following formulation:
\begin{equation}
\begin{split}
    &a(u, v;\mu)+a_{\text{SUPG}}(u, v;\mu)=f(v)+f_{\text{SUPG}}(v),\\
    \text{with }&a_{\text{SUPG}}=\tau_{\text{SUPG}} \int_{\Omega} \left[ -10^{-\mu} (\Delta u + \beta \cdot \nabla u) \cdot (h_{\boldsymbol{x}}\beta \cdot \nabla v) \right]\, d\boldsymbol{x},\\
    \text{and }&f_{\text{SUPG}}= \tau_{\text{SUPG}} \int_{\Omega}  \, (h_{\boldsymbol{x}}\beta \cdot \nabla v)\, d\boldsymbol{x},
\end{split}
    \label{eq:supg}
\end{equation}
where $\tau_{\text{SUPG}}$ is a parameter identifying the weight of the additional term ($\tau_{\text{SUPG}}=0.5$ in our case), while $h_{\boldsymbol{x}}$ is a function defining the diameter associated with each element.

After discretizing the SUPG formulation using a FEM approach, we obtain the following linear system:
\begin{equation}
    (\mathbf{A}_1 + 10^{-\mu}\mathbf{A}_2)\,\mathrm{u}=\mathrm{f}.
    \label{eq:supg-fom-system}
\end{equation}
The full-order quantities are defined as follows:
\begin{itemize}
    \item $(\mathbf{A}_1)_{ij}=\int_{\Omega} [(\beta \cdot \nabla \phi_j)\, \phi_i + \tau_{\text{SUPG}} (\beta \cdot \nabla \phi_j) \cdot (h_{\boldsymbol{x}}\beta \cdot \nabla \phi_i) ]d\boldsymbol{x}$;
    \item $(\mathbf{A}_2)_{ij}=\int_{\Omega} [ \nabla \phi_j \cdot \nabla \phi_i + \tau_{\text{SUPG}} (-\Delta \phi_j) \cdot (h_{\boldsymbol{x}}\beta \cdot \nabla \phi_i) ] d\boldsymbol{x}$;
    \item $(\mathrm{f})_{i}=\int_{\Omega} [\phi_i +  \tau_{\text{SUPG}}  \, (h_{\boldsymbol{x}}\beta \cdot \nabla \phi_i)]\, d\boldsymbol{x}$.
\end{itemize}
The notation $\phi_i$ refers to the FEM basis functions, while the full dimensionality is in this case $N=2500$.

\begin{figure}[htpb!]
    \centering
    \subfloat[$\mathbf{A}_1$]{\includegraphics[width=0.4\linewidth]{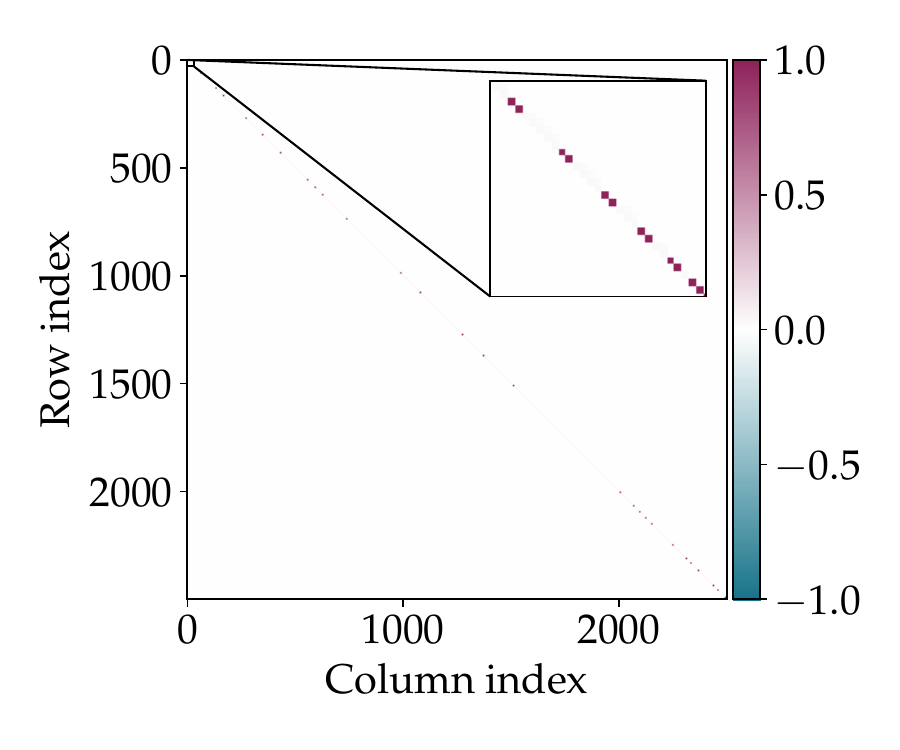}}
    \hfill
    \subfloat[$\mathbf{A}_2$]{\includegraphics[width=0.4\linewidth]{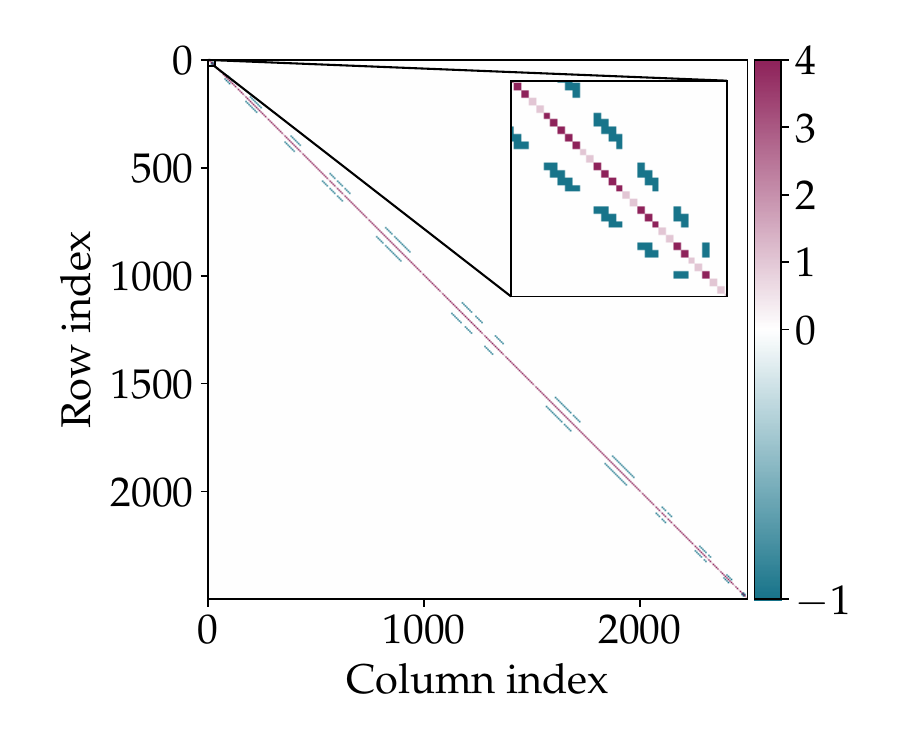}}
    \caption{\emph{Advection-diffusion test case.} Full-order matrices $\mathbf{A}_1$ and $\mathbf{A}_2$.}
    \label{fig:adv-dom-full-matrices}
\end{figure}

The FOM formulation of Equation \eqref{eq:supg-fom-system} allows for a fully separation of the non-parametrized and parametrized matrices $\mathbf{A}_1$ and $\mathbf{A}_2$, which are represented in their sparsity pattern in Figure \ref{fig:adv-dom-full-matrices}.
As in the previous test case, we can recognize patterns in the matrices, e.g., $\mathbf{A}_2$ is symmetric. However, since in the Poisson case the \emph{structure-preserving} approach (S-CE-AE) and the standard CE-AE converged to similar results, we proceed now without imposing any constraint.

The resulting ROM is written as:
\begin{equation}
    (\mathbf{A}_{r1} + 10^{-\mu}\mathbf{A}_{r2})\,\mathbf{z}=\mathbf{y},
    \label{eq:supg-rom-system}
\end{equation}
with $r=2,3,$ or $4$ as latent dimension.

Hence, in this test case we do not need the addition of an interpolation or regression strategy that maps the parameters' space into the reduced quantities.

\begin{figure}[htpb!]
    \centering
    \subfloat[Parameters' space]{\includegraphics[width=0.55\linewidth]{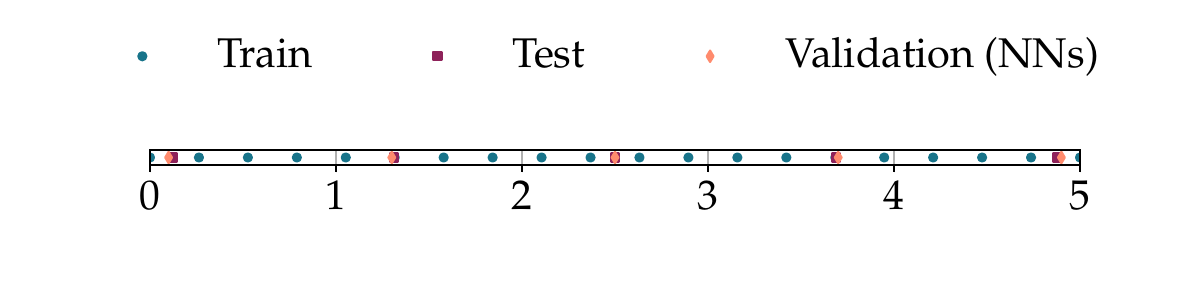}}
    \subfloat[POD singular values' decay]{\includegraphics[width=0.45\linewidth]{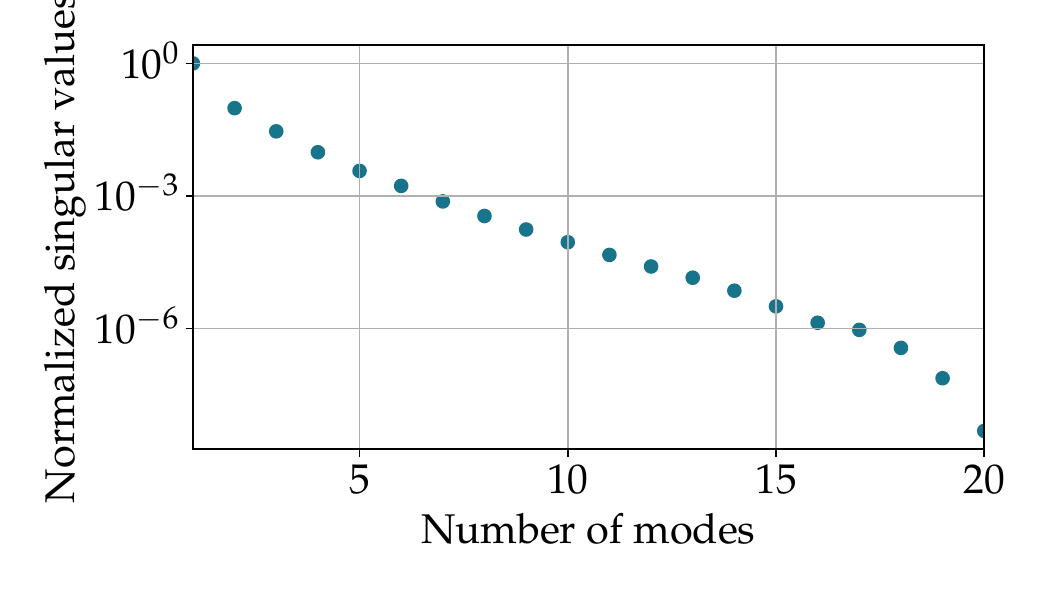}}
    \caption{\emph{Advection-diffusion test case.} Train, validation (for NNs), and test parameters (left), and the POD singular values decay (right).}
    \label{fig:adv-dom-params}
\end{figure}

The parameter values considered for training, validation (only in models involving NNs) and testing, together with the POD singular values' decay, are represented in Figure \ref{fig:adv-dom-params}.
We consider $20$ training parameters, and $5$ parameters both for validation and test. The singular values' decay is similar to the previous test case, highlighting that we need a large number of modes to capture the solution accurately.

\begin{figure}[htpb!]
    \centering
    \subfloat[Latent dimension $r=2$]{\includegraphics[width=0.85\linewidth]{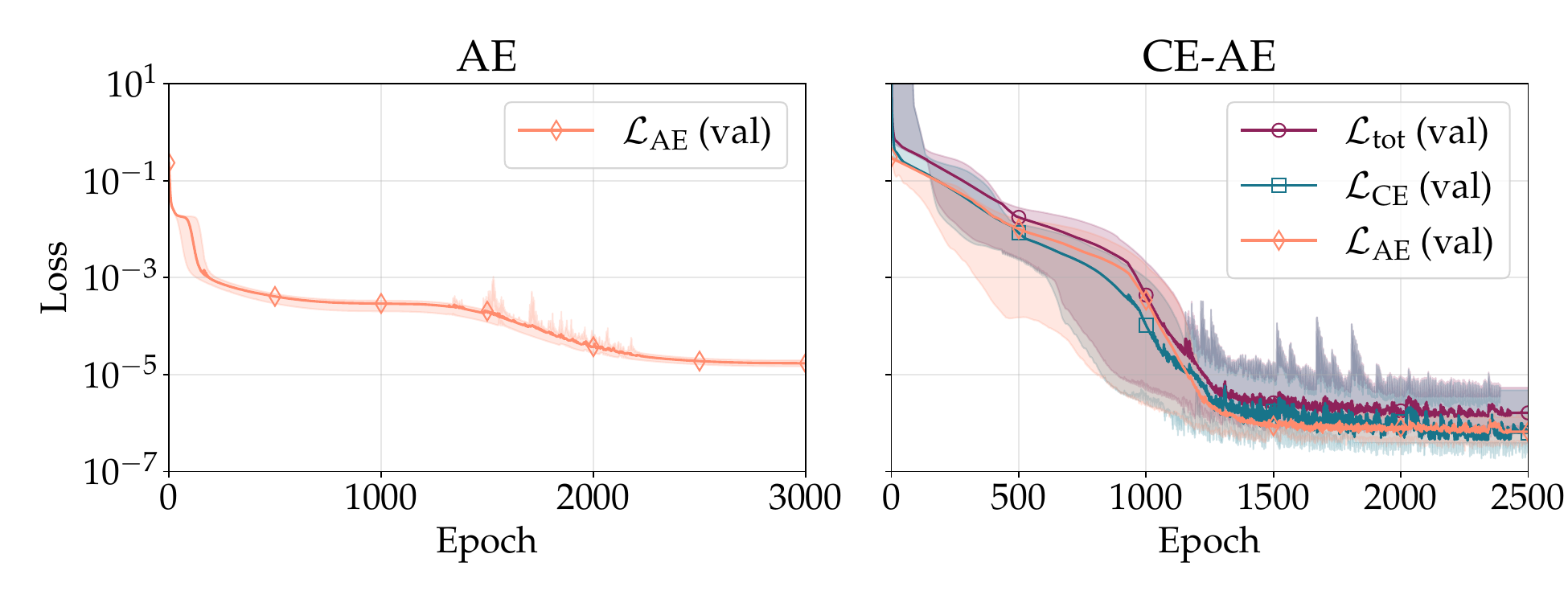}}\\
    \subfloat[Latent dimension $r=3$]{\includegraphics[width=0.85\linewidth]{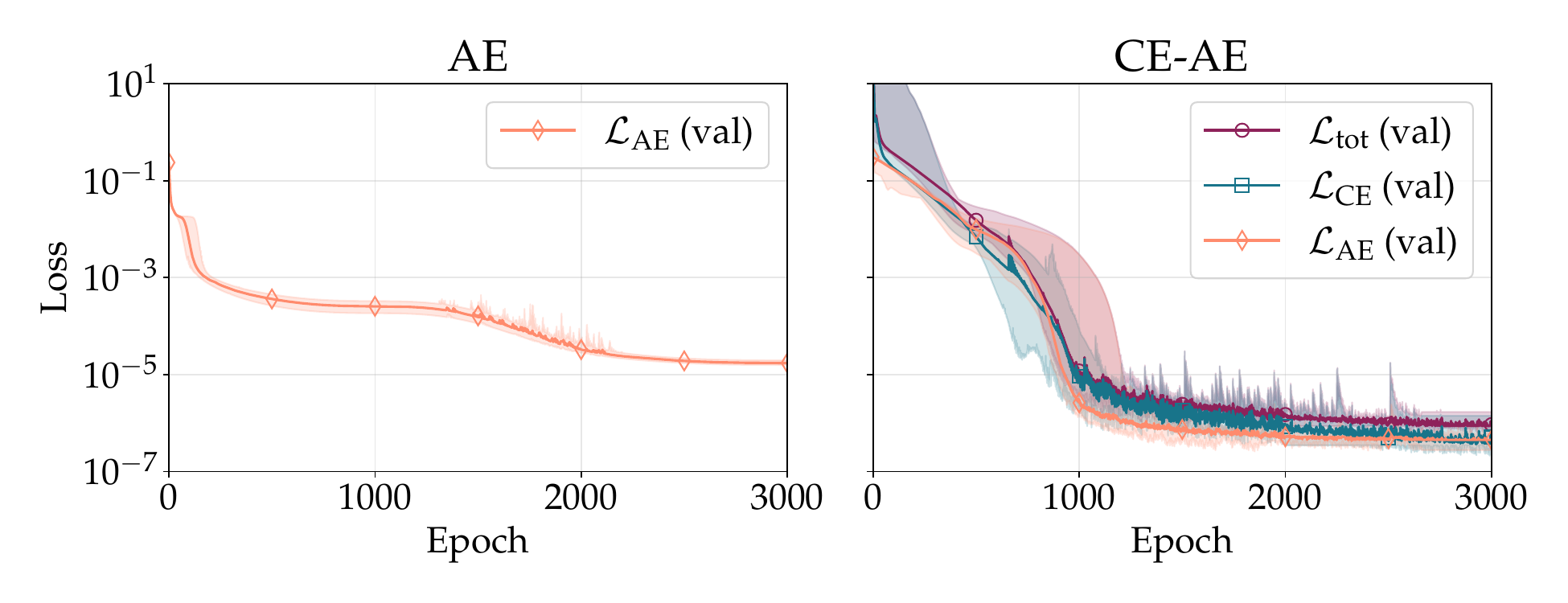}}\\
    \subfloat[Latent dimension $r=4$]{\includegraphics[width=0.85\linewidth]{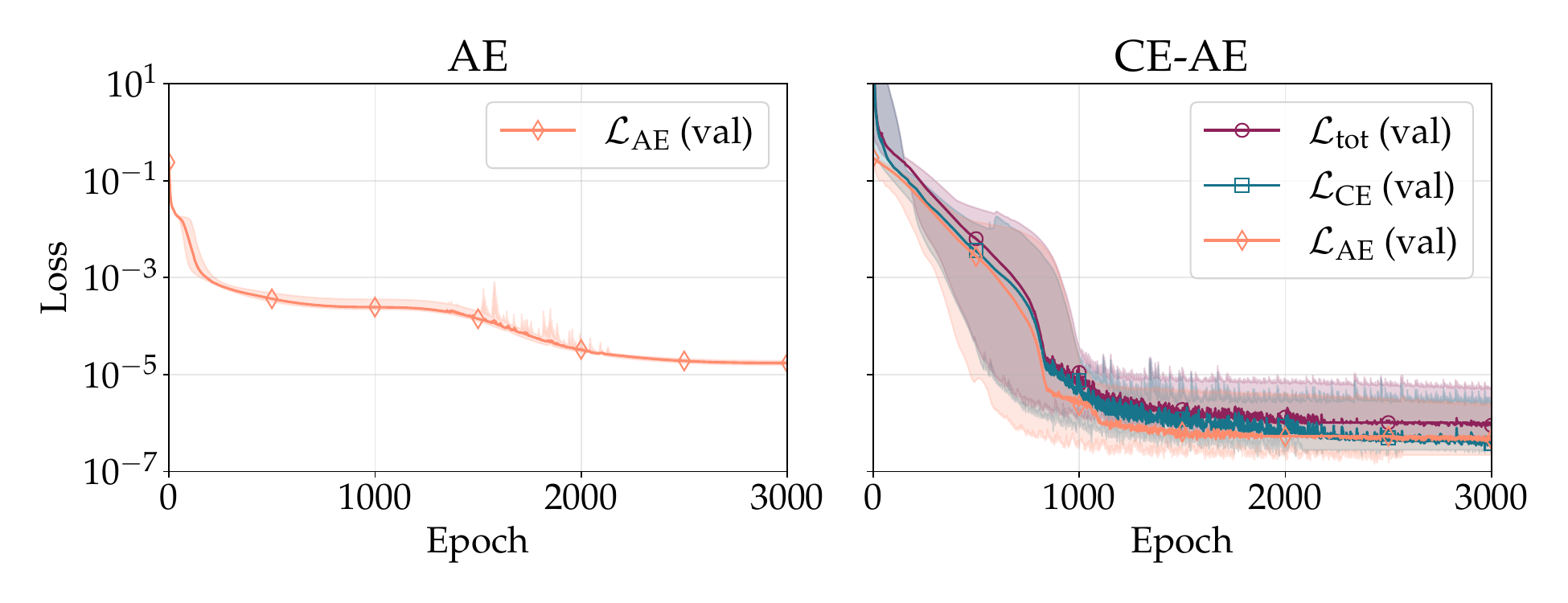}}
    \caption{\emph{Advection-diffusion test case.} Loss trend across epochs for a standard AE, and for the novel model CE-AE, at different latent dimensions, namely $r=2,3,$ and $4$.}
    \label{fig:adv-dom-losses}
\end{figure}

As done in the first test case, we analyse in Figure \ref{fig:adv-dom-losses} the loss behaviour during training iterations, at $r=2,3,4$, for the standard AE and for the CE-AE model.
In all the latent dimensions considered, the CE-AE loss reaches values smaller than more than one order of magnitude with respect to the standard AE. This phenomenon suggests that the convergence is improved by the novel ROM architecture proposed. Moreover, the addition of the loss contribution $\mathcal{L}_{\text{CE}}$ acts as a \emph{stabilizer} for the AE itself, which converges to a more accurate latent representation.

\begin{figure}[htpb!]
    \centering
    \includegraphics[width=\linewidth, trim={3cm 0 5cm 0}, clip]{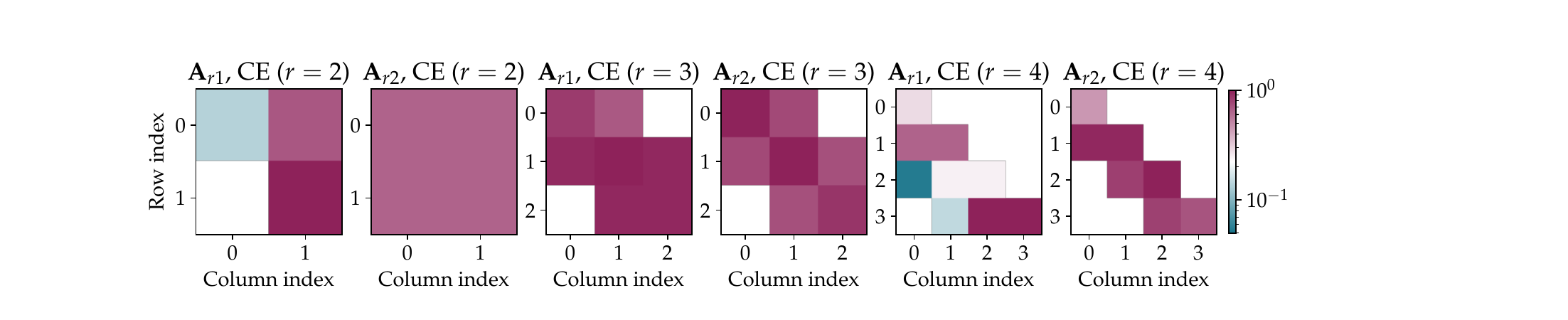}
    \caption{\emph{Advection-diffusion test case.} Reduced matrices $\mathbf{A}_{r1}$ and $\mathbf{A}_{r2}$, normalized in range $[0, 1]$. The matrices are obtained by compressing the full-order discrete matrix $\mathbf{A}_1$ and $\mathbf{A}_2$ with the convolutional encoder CE, at different latent dimensions $r=2, 3,$ and $4$.}
    \label{fig:adv-dom-matrices}
\end{figure}

The reduced matrices obtained by CE compression, and then normalized in the range $[0, 1]$ are represented in Figure \ref{fig:adv-dom-matrices}.

\begin{figure}[htpb!]
    \centering
    \includegraphics[width=\linewidth]{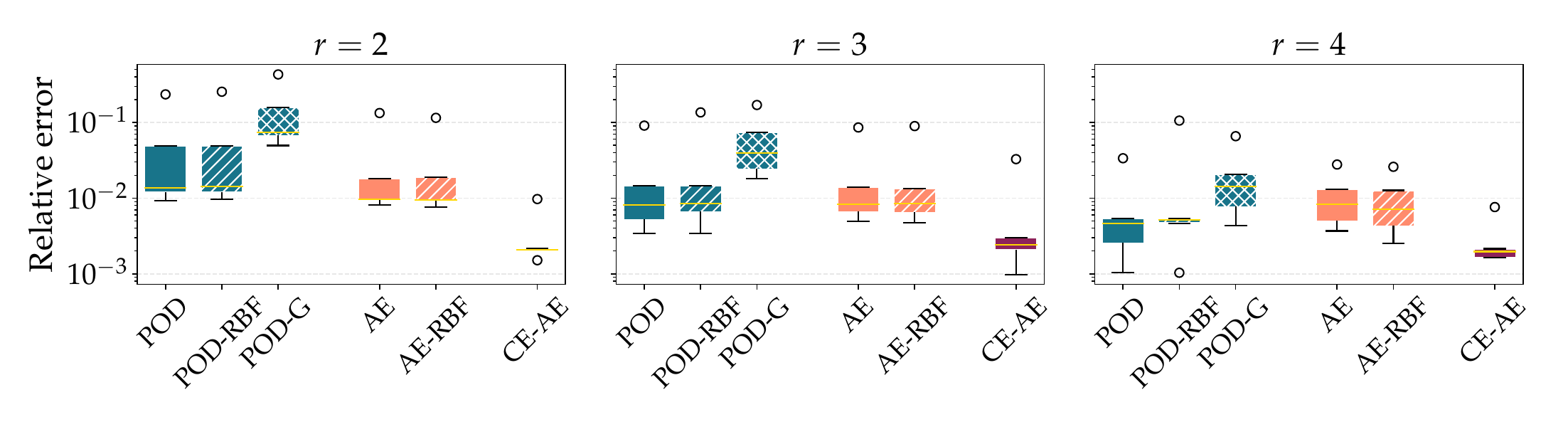}
    \caption{\emph{Advection-diffusion test case.} Relative $L^2$ test errors of different approximations with respect to the FOM reference, at the latent dimensions $r=2, 3$, and $4$.}
    \label{fig:adv-dom-boxplot}
\end{figure}

\begin{figure}[htpb!]
    \centering
    \includegraphics[width=\linewidth, trim={0 2cm 0 1cm}, clip]{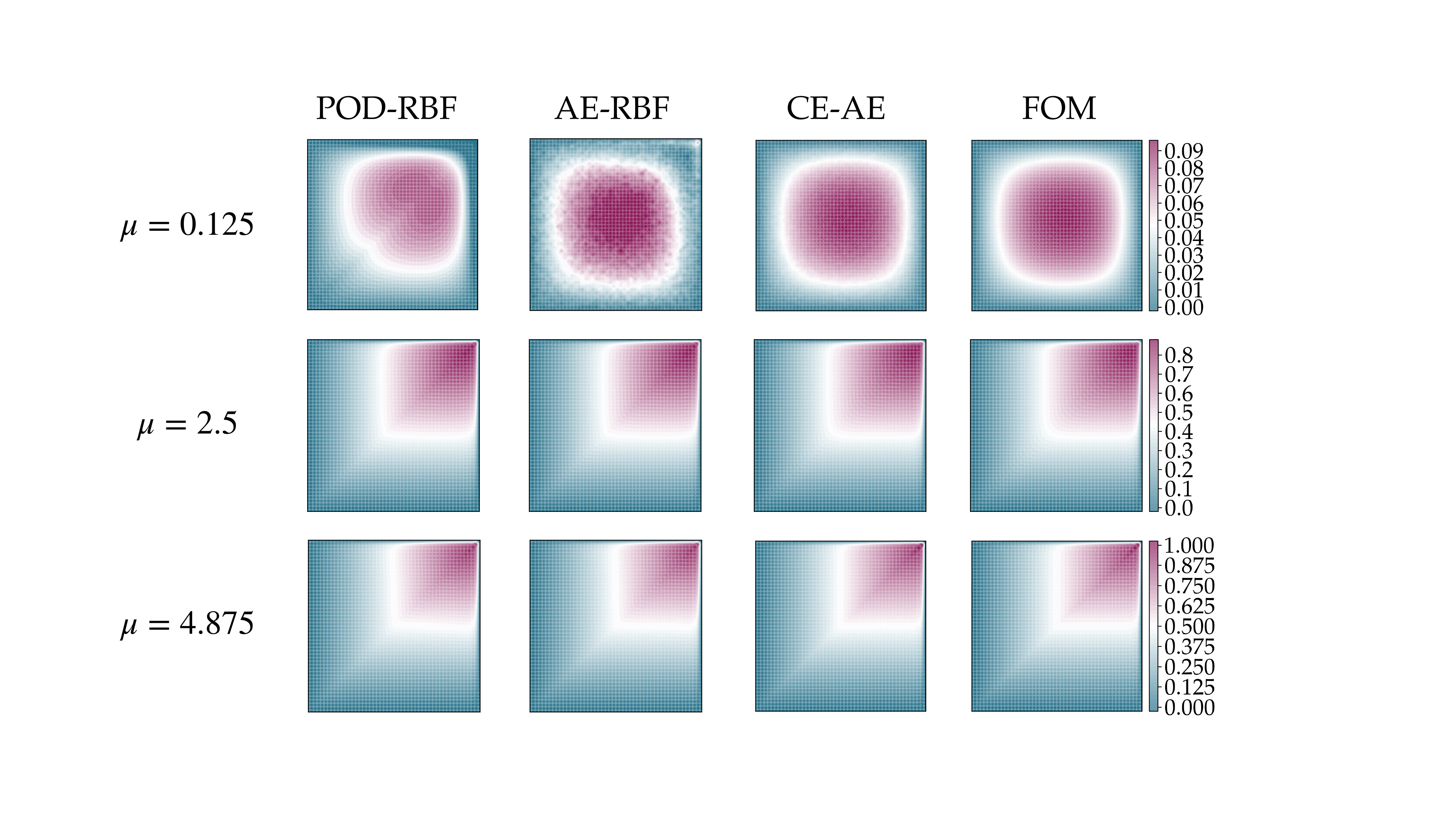}
    \caption{\emph{Advection-dominated test case.} Approximated solutions for POD-RBF, AE-RBF, and CE-AE, compared with the FOM reference, at latent dimension $r=2$, and considering three test parameters.}
    \label{fig:adv-dom-fields}
\end{figure}

At both $r=2$ and $r=3$, the CE converges to symmetric matrices $\mathbf{A}_{r2}$, reproducing the FOM properties. This is not equally satisfied at $r=4$, where the increased latent dimension may lead to more variability in the converged solution.
Moreover, at $r=2$ all the components of $\mathbf{A}_{r2}$ have the same value, suggesting that it is treated as one additional degree of freedom, properly affecting the effect of the parameters on the latent solution.
The POD, AE and CE-AE solutions are quantitatively compared in terms of the relative test $L^2$ error in Figure \ref{fig:adv-dom-boxplot}.
The corresponding fields are graphically compared in Figure \ref{fig:adv-dom-fields} for the POD-RBF, AE-RBF, and CE-AE methods. The results in both figures have been obtained considering the most accurate AE and CE-AE from the ensemble considered to generate Figure \ref{fig:adv-dom-losses}.
At all latent dimensions, the CE-AE is significantly more accurate than the standard POD and AE, and it is characterized by a smaller variability of the error among the test parameters considered. Interestingly, the AE has similar accuracy to the POD, while it is even less accurate at $r=4$.

The accuracy improvement can be also noticed qualitatively in Figure \ref{fig:adv-dom-fields}, which shows the approximated solutions at $r=2$. Here, we can see that the AE and the POD fail to accurately approximate the solution at small $\mu$ values (first row of the figure), while the CE-AE is able to reproduce the FOM field similar to larger $\mu$ values.

\subsection{Burgers test case}
\label{subsec:burgers-cae}
\begin{figure}
    \centering
    \includegraphics[width=0.75\linewidth, trim={0cm 9cm 0 9cm}, clip]{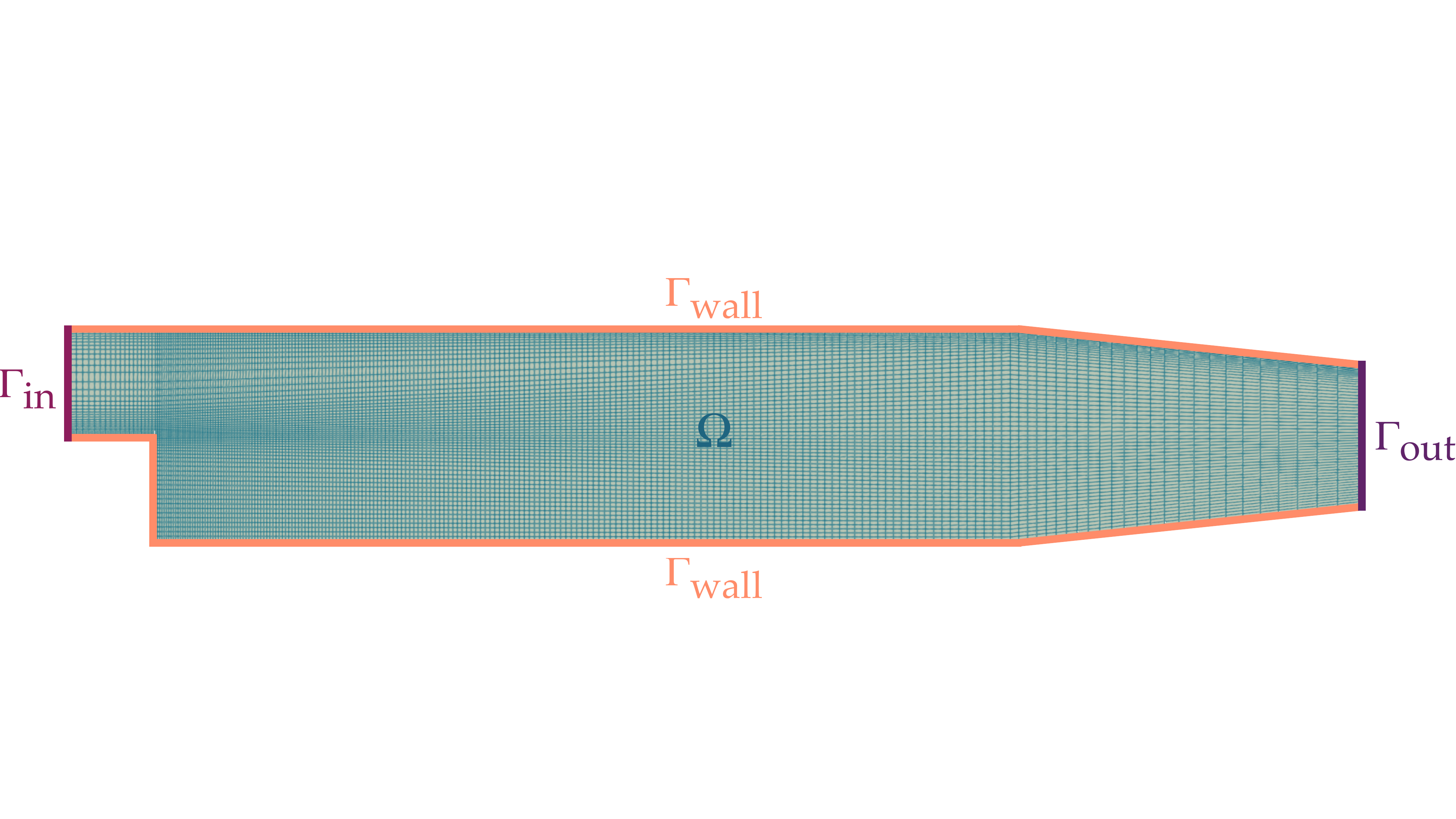}
    \caption{\emph{Burgers test case}. Domain and computational grid adopted, with boundaries' notation.}
    \label{fig:domain-burgers}
\end{figure}

The final test case is a viscous Burgers equation in a 2D backstep domain modeled in OpenFOAM~\cite{jasak1996error, jasak2009openfoam} using the FVM.
The domain and FV grid are represented in Figure \ref{fig:domain-burgers}, including the notation adopted for the boundaries.

The problem in continuous form can be written as follows:
\begin{equation}
\begin{cases}
    \dfrac{\partial \boldsymbol{u}}{\partial t} + (\boldsymbol{u} \cdot \nabla) \boldsymbol{u} - \nu \Delta \boldsymbol{u}=\boldsymbol{0} &\text{ in }\Omega,\\
    \boldsymbol{u}=(1, 0) &\text{ on }\Gamma_{\text{in}},\\
    \nabla \boldsymbol{u}\cdot \boldsymbol{n}=0 &\text{ on } \Gamma_{\text{out}},\\
    \boldsymbol{u}=\boldsymbol{0} &\text{ on }\Gamma_{\text{wall}}.
\end{cases}
    \label{eq:burgers-continuous}
\end{equation}

The continuous problem may be discretized as follows:
\begin{equation}
    \mathbf{M_u}\dfrac{\udisc_{n+1}-\udisc_n}{\Delta t} + \mathbf{C_u}(\udisc_n)\udisc_{n+1} - \nu \mathbf{D_u}\udisc_{n+1} = \mathbf{f_C}(\udisc_n) + \mathbf{f_D},
    \label{eq:burgers-fom}
\end{equation}

where $\mathbf{M_u}$ is the mass matrix, $\mathbf{C_u}(\udisc)$ is the discrete convection operator, $\mathbf{D_u}$ is the laplacian operator, while $\mathbf{f_D}$ and $\mathbf{f_C}$ are the known vectors associated with the laplacian and convection terms, respectively. In practice, the finite-volume discretization of each operator gives rise not only to a sparse matrix but also to a source vector, which incorporates the boundary conditions.
In a more compact form, we can rewrite Equation \eqref{eq:burgers-fom} as follows:
\begin{equation}
    \boxed{\left( \mathbf{M_u} + \Delta t (\mathbf{C_u}(\udisc_n) + \mathbf{D_u})\right)\, \udisc_{n+1} = \mathbf{f_C}+\mathbf{f_D}+\mathbf{M_u}\udisc_n.
    }
    \label{eq:burgers-fom-compact}
\end{equation}

Differently to the previous cases, the FOM number of degrees of freedom is now larger ($N=24450$, with $12250$ cells), making the memory demand to compress the operators $\mathbf{M_u}$, ${\mathbf{C_u}}$ and $\mathbf{D_u} \in \mathbb{R}^{N\times N}$ prohibitively high.
For this reason, we compress the FOM matrices using a continuous convolutional encoder to exploit the matrices sparsity, as introduced in Section \ref{subsec:cceae}.

The resulting ROM system can be written as follows:
\begin{equation}
   \boxed{\left( \mathbf{M_u}_r + \Delta t (\mathbf{C_u}_r(\mathbf{z}_n) + \mathbf{D_u}_r)\right)\, \mathbf{z}_{n+1} = \mathbf{y_C}(\mathbf{z}_n)+\mathbf{y_D}+\mathbf{M_u}_r{\mathbf{z}}_n.
    }
    \label{eq:burgers-rom-system}
\end{equation}
To clarify how the reduced system is constructed, we describe the pseudo-code of the offline stage in Algorithm \ref{alg-burgers-rom-offline} and the online prediction in Algorithm \ref{alg-burgers-rom-online}. In the \textbf{offline} stage, the logic adopted is the same as in the previous test cases, with some differences:
\begin{itemize}
    \item There are multiple matrices to compress, and we adopt the same convolutional encoder. Other choices are possible (different channels, multiple encoders) but we postpone this investigation to future work.
    \item We only adopt an intrusive loss $\mathcal{L}_{\text{CCE}}$ to train the CCE-AE model. In our preliminary investigations, we noticed that including both the losses contributions led to slow or no convergence during training. Therefore, we decided to proceed without including the purely data-driven contribution into the loss.
    \item We deal with a time-dependent problem with fixed $\Delta t=\num{1e-4}$, which is kept unchanged in the ROM formulation (Equation \eqref{eq:burgers-rom-system}).
    \item The RBF interpolation (second part of Algorithm \ref{alg-burgers-rom-offline}) is applied to the time-dependent matrix $\mathbf{C_u}_r$ and vector $\mathbf{y_C}_r$. Since in the FOM such quantities are functions of the FOM state $\udisc_n$, we reproduce the same maps at the reduced level, namely $\mathbf{C_u}_r(\mathbf{z}_n)$ and $\mathbf{y_C}_r(\mathbf{z}_n)$. Therefore, we distinguish two RBF interpolators, one for the matrices (RBF$_{\text{mat}}$), which interpolates the matrices element-wise, and one for the vector (RBF$_{\text{vec}}$), which take as input the latent state rather than the parameter (e.g., time). This allows the construction of more interpretable latent maps, as we will see in the results.
    \item The input used to train the RBF maps is composed of the latent trajectory $\mathbf{z}^{\text{train}}_n$, $n=1, N_{\mu}$, where the $\mathbf{z}^{\text{train}}$ have been found by compressing the FOM quantities with the pre-trained CCE and $\mathcal{E}$, and solving the reduced system.
\end{itemize}

In the \textbf{online} stage, we compress the time-invariant matrices $\mathbf{M_u}_r$ and $\mathbf{D_u}_r$ using the CCE, while the time-dependent $\mathbf{C_u}_r(t_n)$ are predicted at each time step $t_n$ using the pre-trained RBF$_{\text{mat}}$ map having as input the reduced solution $\mathbf{z}_{n}$. The same applies to the reduced vectors: $\mathbf{y_D}$ is found once by compressing the time-invariant vector $\mathbf{f_D}$, while the time-dependent vector $\mathbf{y_C}(t_n)$ is found as RBF$_{\text{vec}}(\mathbf{z}_n)$ at each time step.

\begin{algorithm*}[htpb!]
\caption{Pseudo-code for the offline training of the CCE-AE-RBF approach.}
\label{alg-burgers-rom-offline}
\begin{algorithmic}[1]
\Statex
\Statex \textbf{\textcolor{red}{Phase I: CCE-AE training}}
\State \textbf{Input data}: snapshots $\udisc_n$, $n=1, \dots, N_\mu$; matrices $\mathbf{M_u}$, $\mathbf{D_u}$, and $\mathbf{C_u}(\udisc_n)$, $n=1, \dots, N_{\mu}$; vectors $\mathbf{f_D}$, $\mathbf{f_C}(\udisc_n)$, $n=1, \dots, N_{\mu}$;
\State Continuous convolutional encoder CCE; MLP encoder $\mathcal{E}$ and decoder $\mathcal{D}$; number of training epochs $n_{\text{epochs}}$;
\For{$e \in [1, \dots, n_{\text{epochs}}]$}
\State \textbf{Initial FOM state}:
\State $\mathbf{M_u}_r= \text{CCE}(\mathbf{M_u})$;
\State $\mathbf{D_u}_r=\text{CCE}(\mathbf{D_u})$;
\State $\mathbf{y_D}=\mathcal{E}(\mathbf{f_D})$;
\State $\mathbf{z}_0 = \mathcal{E}(\mathbf{u}_0)$;
\For{$n \in [0, \dots, N_{\mu} - 1]$}
\State $\mathbf{C_u}_r \gets \text{CCE}(\mathbf{C_u}(\udisc_n))$;
\State $\mathbf{y_C} \gets \mathcal{E}(\mathbf{f_C}(\udisc_n))$;
\State $\mathbf{A}_r \gets \mathbf{M_u}_r + \Delta t (\mathbf{C_u}_r + \mathbf{D_u}_r)$;
\State $\mathbf{y} \gets \mathbf{y_C}+\mathbf{y_D}+\mathbf{M_u}_r{\mathbf{z}}_n$;
\State Compute $\mathbf{z}_{n+1}$ as the solution of $\mathbf{A}_r \mathbf{z}_{n+1}=\mathbf{y}$;
\EndFor
\State Compute the loss $\mathcal{L}_{\text{CCE}}=\dfrac{1}{N_{\mu}}\sum_{n=1}^{N_{\mu}} \|\mathcal{D}(\mathbf{z}_{n})-\mathbf{u}_n \|_2^2$.
\State Backpropagate $\mathcal{L}_{\text{CCE}}$ and update the trainable parameters of CCE, $\mathcal{E}$, and $\mathcal{D}$;
\EndFor
\Statex
\Statex \textbf{\textcolor{blue}{Phase II: RBF training}}
\State \textbf{Collect input}: pre-trained $\mathcal{E}$ and CCE;
\State $\mathbf{M_u}_r= \text{CCE}(\mathbf{M_u})$;
\State $\mathbf{D_u}_r=\text{CCE}(\mathbf{D_u})$;
\State $\mathbf{y_D}=\mathcal{E}(\mathbf{f_D})$;
\State $\mathbf{C_u}^{\text{train}}_r(t_n) = \text{CCE}(\mathbf{C_u}(\udisc_n))$ for $n=0, \dots, N_{\mu}$;
\State $\mathbf{y_C}^{\text{train}}(t_n) = \mathcal{E}(\mathbf{f_C}(\udisc_n))$ for $n=0, \dots, N_{\mu}$;
\State $\mathbf{z}^{\text{train}}_0 = \mathcal{E}(\mathbf{u}_0)$;
\For{$n \in [0, \dots, N_{\mu} - 1]$}
\State $\mathbf{A}_r \gets \mathbf{M_u}_r + \Delta t (\mathbf{C_u}_r + \mathbf{D_u}_r)$;
\State $\mathbf{y} \gets \mathbf{y_C}+\mathbf{y_D}+\mathbf{M_u}_r{\mathbf{z}}^{\text{train}}_n$;
\State Compute $\mathbf{z}^{\text{train}}_{n+1}$ as the solution of $\mathbf{A}_r \mathbf{z}^{\text{train}}_{n+1}=\mathbf{y}$;
\EndFor
\State \textbf{Input (RBF train)}: $[\mathbf{z}^{\text{train}}_{1}, \dots, \mathbf{z}^{\text{train}}_{N_{\mu}}]$, $[\mathbf{y_C}^{\text{train}}(t_1), \dots, \mathbf{y_C}^{\text{train}}(t_{N_{\mu}})]$, $[\mathbf{C_u}^{\text{train}}_r(t_1), \dots, \mathbf{C_u}^{\text{train}}_r(t_{N_{\mu}})]$;
\State Train RBF$_{\text{mat}}: \mathbf{z}^{\text{train}}_n \mapsto \mathbf{C_u}^{\text{train}}_r(t_n)$;
\State Train RBF$_{\text{vec}}: \mathbf{z}^{\text{train}}_n \mapsto \mathbf{y_C}^{\text{train}}(t_n)$.
\end{algorithmic}
\end{algorithm*}

\begin{algorithm*}[htpb!]
\caption{Pseudo-code for the online prediction of the CCE-AE-RBF approach.}
\label{alg-burgers-rom-online}
\begin{algorithmic}[1]
\State \textbf{Input}: Pre-trained continuous convolutional encoder CCE, MLP encoder $\mathcal{E}$ and decoder $\mathcal{D}$;
\State \textbf{Initial FOM state}: $\mathbf{u}_0=\mathbf{u}(t_0)$, time steps $t_n$, $n=1, \dots N_T$
\State $\mathbf{M_u}_r= \text{CCE}(\mathbf{M_u})$;
\State $\mathbf{D_u}_r=\text{CCE}(\mathbf{D_u})$;
\State $\mathbf{y_D}_r=\mathcal{E}(\mathbf{f_D})$;
\State $\mathbf{z}_0 = \mathcal{E}(\mathbf{u}_0)$;
\For{$n \in [0, \dots, N_T - 1]$}
\State $\mathbf{C_u}_r \gets \text{RBF}_{\text{mat}}(\mathbf{z}_n)$;
\State $\mathbf{y_C}_r \gets \text{RBF}_{\text{vec}}(\mathbf{z}_n)$;
\State $\mathbf{A}_r \gets \mathbf{M_u}_r + \Delta t (\mathbf{C_u}_r + \mathbf{D_u}_r)$;
\State $\mathbf{y} \gets \mathbf{y_C}+\mathbf{y_D}+\mathbf{M_u}_r{\mathbf{z}}_n$;
\State Compute $\mathbf{z}_{n+1}$ as the solution of $\mathbf{A}_r \mathbf{z}_{n+1}=\mathbf{y}$.
\State Reconstruct the solution as $\mathcal{D}(\mathbf{z}_{n+1})$.
\EndFor
\end{algorithmic}
\end{algorithm*}

\begin{figure}[htpb!]
    \centering
    \includegraphics[width=0.55\linewidth]{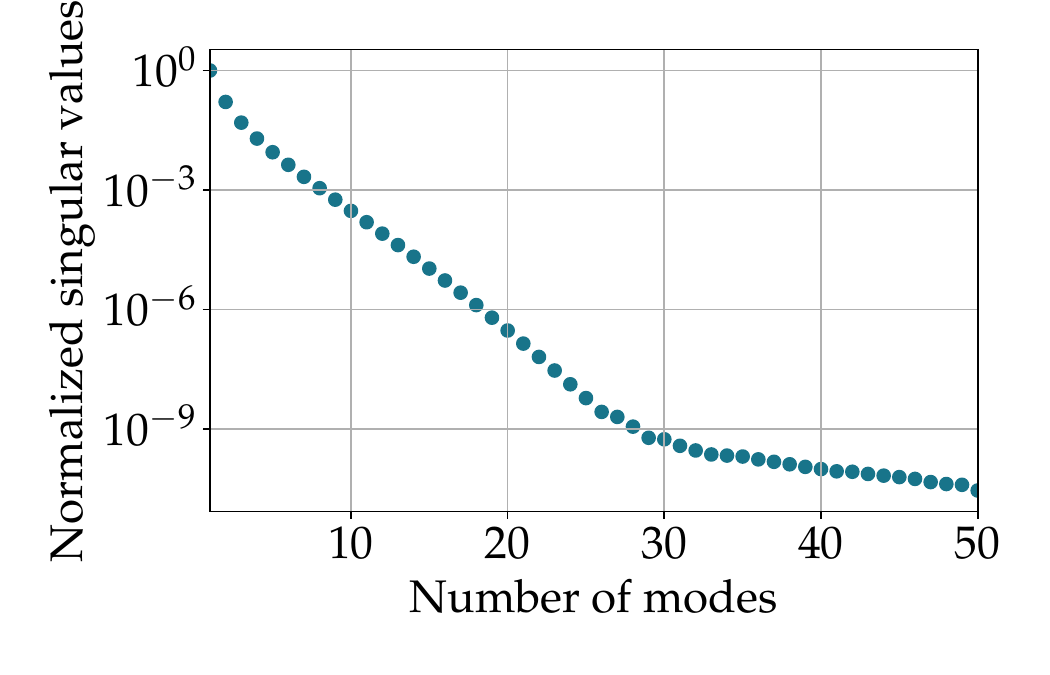}
    \caption{\emph{Burgers test case}. The POD singular values decay.}
    \label{fig:sing-vals-burgers}
\end{figure}

As previously mentioned, the only parameter is time, considering the first $N_{\mu}=50$ training time steps, and testing on the window $[0, N_T]$, with $N_T=100$. The time steps $t_n$, with $n \in [50, 75]$ are used during training for validation.
If we apply the POD to the snapshots containing the first $N_{\mu}=50$ snapshots, we obtain the singular values' decay represented in Figure \ref{fig:sing-vals-burgers}.

\begin{figure}[htpb!]
    \centering
    \subfloat[Latent dimension $r=2$]{\includegraphics[width=0.9\linewidth]{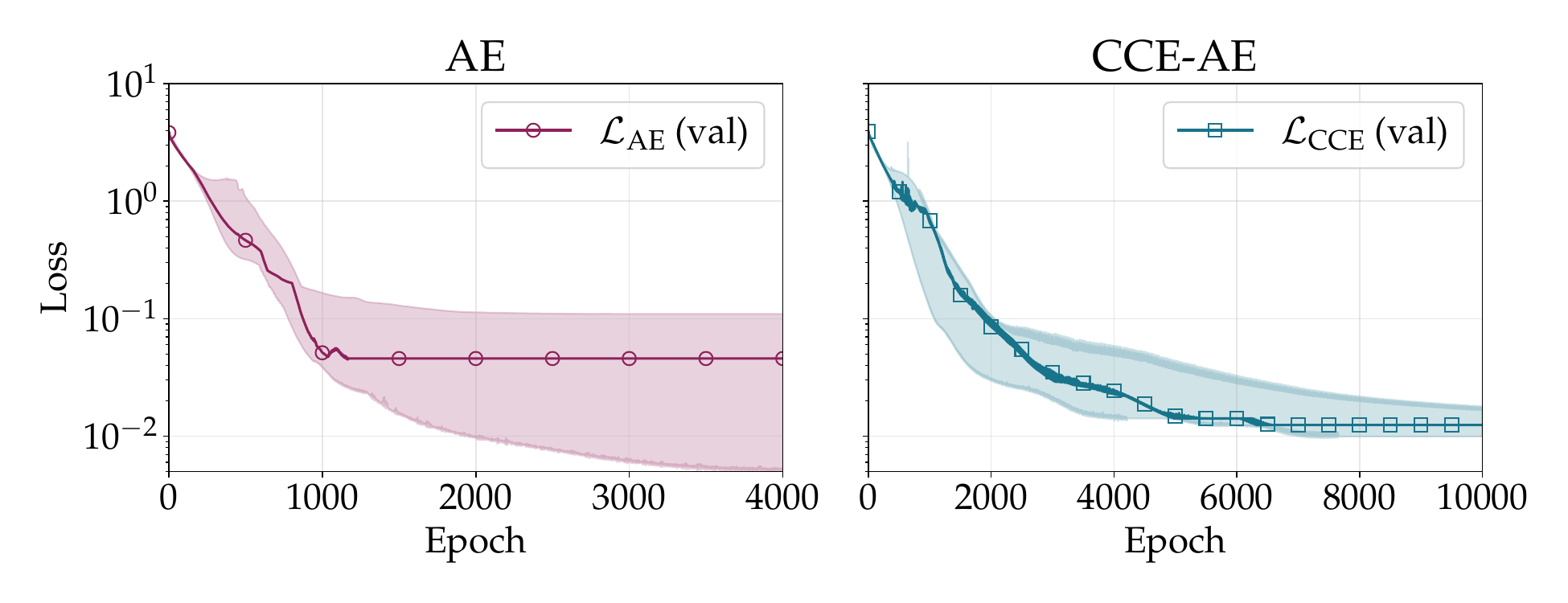}}\\
    \subfloat[Latent dimension $r=3$]{\includegraphics[width=0.9\linewidth]{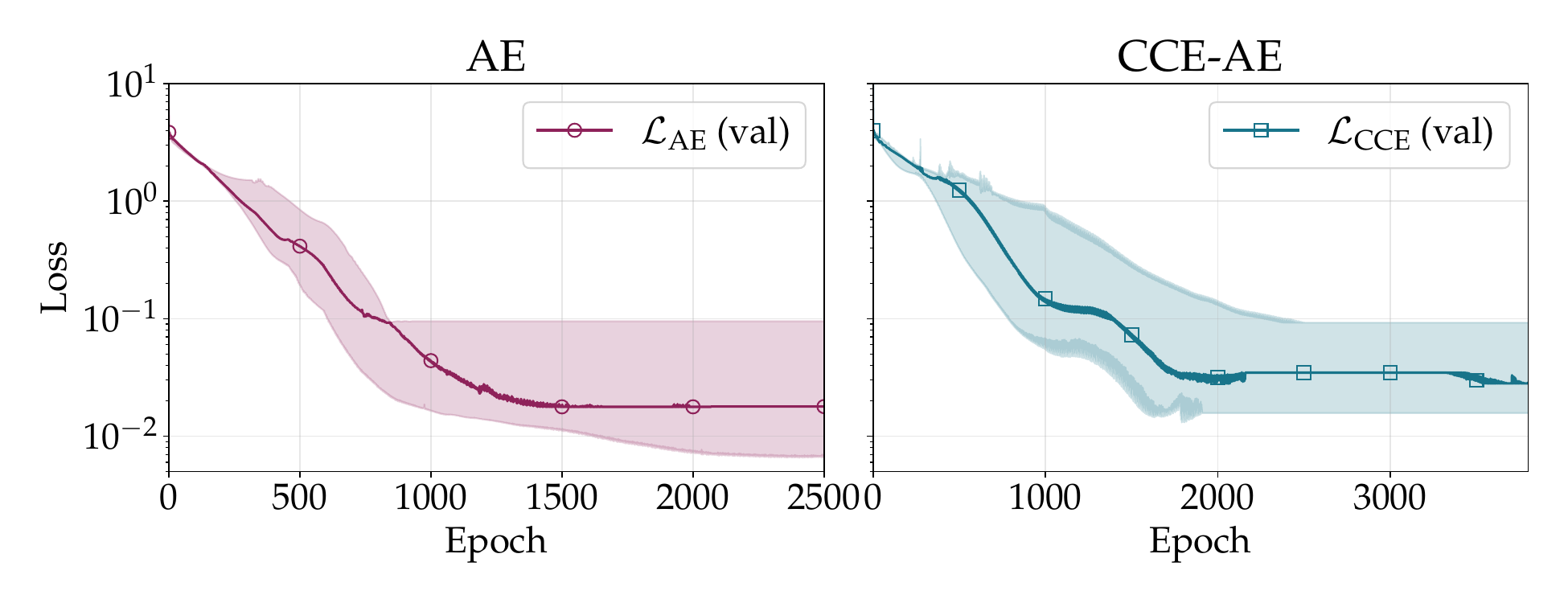}}
    \caption{\emph{Burgers test case.} Loss trend across epochs for a standard AE, and for the novel model CCE-AE, at different latent dimensions, namely $r=2$ and $3$.}
    \label{fig:burgers-losses}
\end{figure}

Figure \ref{fig:burgers-losses} shows the validation loss decay during training, comparing the training process of a standard AE and of the proposed CCE-AE approach.
Since this test case is particularly challenging, both the AE and CCE-AE validation losses converge to relatively high values (approximately $10^{-2}$). This indicates that the extrapolation task is inherently difficult, as even the unconstrained AE struggles to accurately reconstruct unseen samples. However, we can observe that, also in this case, for the smallest latent dimension ($r=2$), the CCE-AE validation loss converges to a lower value than that of the AE. This suggests that coupling the autoencoder with the consistency constraint improves the robustness of the training process.

\begin{figure}[htpb!]
    \centering
    \includegraphics[width=\linewidth]{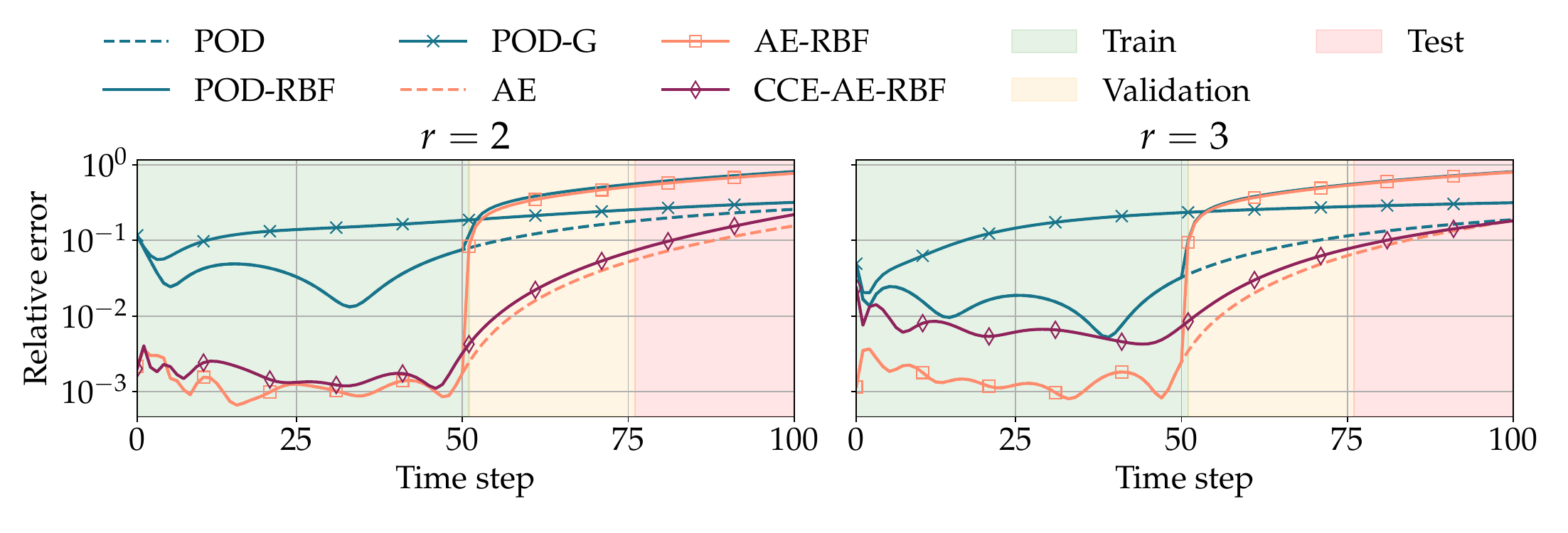}
    \caption{\emph{Burgers test case}. Time variation of the relative $L^2$ errors with respect to the FOM solution, obtained using standard methodologies (POD, POD-RBF, POD-G, AE, AE-RBF), and the novel CCE-AE-RBF approach.}
    \label{fig:burgers-time-errs}
\end{figure}

We now analyse the predictive accuracy of the proposed model, by testing it in a time extrapolation setting, namely at $t_n, n \in [75, 100]$. Figure \ref{fig:burgers-time-errs} represents the time evolution of the relative $L^2$ errors of the approximated solution with respect to the FOM solution, for different POD-based, and AE-based approaches.
As can be seen, linear approaches such as the POD fail to accurately reconstruct the solution in an under-resolved regime. Introducing nonlinearities significantly improves the reconstruction capability: within the training window, the AE achieves reconstruction errors one order of magnitude smaller than those of the POD for both $r=2$ and $r=3$. Nevertheless, when the solution evolves beyond the training regime, the AE fails to preserve the correct temporal dynamics, resulting in inaccurate extrapolations.
Introducing reduced-order interpolation through RBF (POD-RBF and AE-RBF) further deteriorates the reconstruction accuracy. This can be attributed to the difficulty of approximating the mapping from the temporal variable $t$ to the reduced coefficients $\mathbf{z}(t)$. In fact, the reduced trajectories may exhibit complex and nonlinear behaviours that cannot be accurately represented by the chosen interpolation approach. Consequently, these approaches become unreliable for prediction beyond the training regime and cannot be effectively employed for extrapolation.

A POD-Galerkin approach is also analysed. We adopt an \emph{inefficient} formulation, in which the full-order quantities $\mathbf{C_u}$ and $\mathbf{f_C}$ are projected onto the POD modes, also in the predictive regime. Despite the use of full-order quantities, the method already exhibits large errors from the beginning of the prediction window, indicating that the reduced formulation is not able to accurately capture the solution dynamics.

Finally, the CCE-AE-RBF approach is considered. Although the error still increases over time beyond the training window, the proposed interpolation strategies $\mathbf{z}(t)\mapsto \mathbf{y_C}(t)$ and $\mathbf{z}(t)\mapsto \mathbf{C_u}_r(t)$ provide a more accurate approximation of the reduced quantities, keeping the prediction error closer to the reconstruction error. This suggests that the main limitation is not introduced by the RBF interpolation itself, but rather by the difficulty of the autoencoder in representing and propagating the nonlinear dynamics of the solution.

We can summarize the main findings here:
\begin{itemize}
    \item {The main source of error appears to be related to the autoencoder representation.} Due to the strongly nonlinear nature of the Burgers dynamics, the latent space learned by the AE is not sufficiently accurate to preserve the temporal evolution of the solution in the predictive regime. Therefore, the observed degradation is mainly attributed to the difficulty of the underlying nonlinear reduction rather than to the proposed architecture.
\item {The RBF interpolation strategy provides a more suitable approximation framework.} Unlike approaches based on the direct mapping $t \mapsto \mathbf{z}(t)$, the proposed method interpolates the reduced operators as a function of the latent coordinates. This mapping is more consistent with the reduced-order formulation and introduces a smaller additional error, allowing the predictions to remain closer to the intrinsic reconstruction capability of the CCE-AE.
\end{itemize}

\begin{figure}[htpb!]
    \centering
    \subfloat[$x$ velocity component]{\includegraphics[width=\linewidth, trim={5cm 20cm 5cm 5cm}, clip]{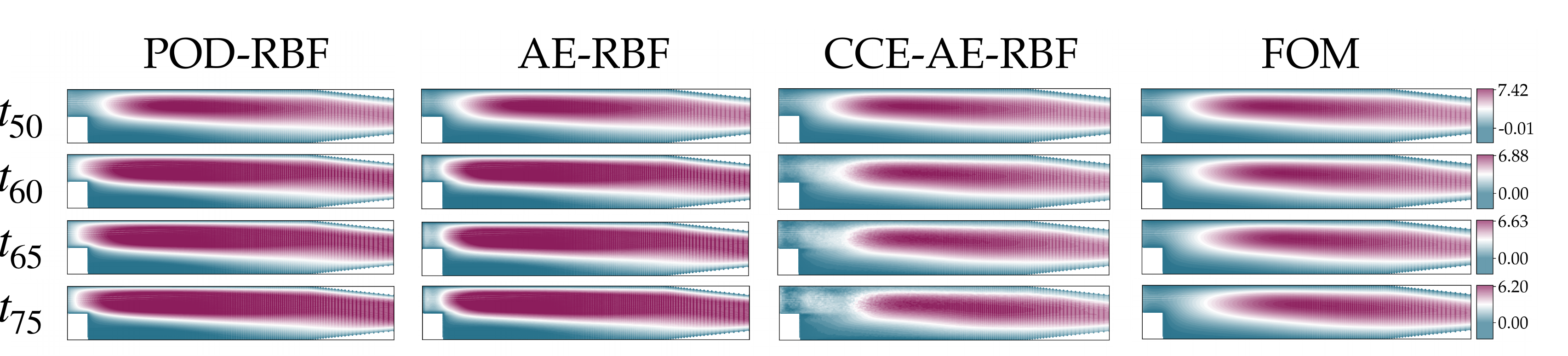}}\\
    \subfloat[$y$ velocity component]{\includegraphics[width=\linewidth, trim={5cm 20cm 5cm 5cm}, clip]{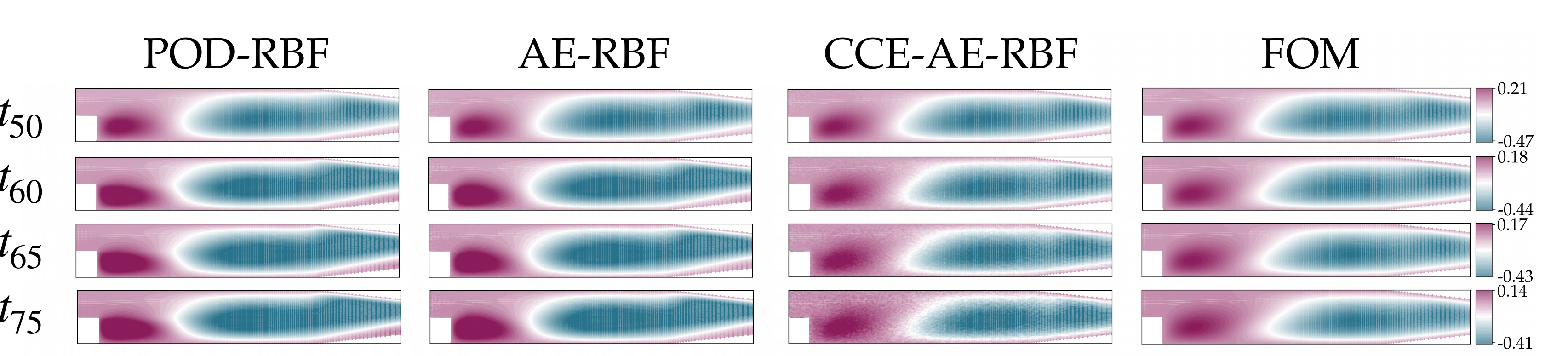}}
    \caption{\emph{Burgers test case}. Solution fields for the $x$ (top) and $y$ (bottom) components, found using POD-RBF, AE-RBF, CCE-AE-RBF, compared with the FOM. The fields are shown at different time steps in a predictive regime, at the fixed latent dimension $r=3$.}
    \label{fig:fields-burgers}
\end{figure}

Finally, Figure \ref{fig:fields-burgers} compares the graphical approximations obtained by using different ROMs at fixed latent dimension $r=3$, at different time steps (in the predictive setting).
As expected, the POD-RBF and AE-RBF approximations are totally unphysical, and are not following the time evolution of the FOM solution. On the other hand, the CCE-AE-RBF approach, although characterized by some spurious oscillations, is able to correctly reproduce the overall behaviour of the FOM, suggesting that the physics-based ROM formulation provides a more accurate representation of the underlying dynamics. The oscillations may be related to the limited capacity of the decoder, which may require a larger architecture (e.g., more layers and neurons). Different decoder architectures will be explored in future works.

\subsection{Computational details}
\label{subsec:cpu-time}
This section describes the computational aspects of the proposed methods, including the architecture and hyperparameters' setup of the neural networks, and the training cost (number of parameters and wall time).

The training process has been performed using the SOAP optimizer~\cite{vyas2025soap} using one GPU NVIDIA RTX A2000 (12 GB VRAM), CUDA 12.4.

The AE is a deep fully-connected neural network with input and output dimensions equal to $N$. For the first two test cases, the encoder consists of hidden layers with $[100, 30, r]$ neurons, whereas for the Burgers test case we adopt a more compact architecture, $[10, r]$. In both cases, the decoder has a symmetric architecture to the encoder, and we adopt the Softplus as activation function.
Preliminary numerical experiments showed that increasing the network size for the Burgers problem does not provide any significant improvement in reconstruction accuracy, while it increases the risk of overfitting. Therefore, a smaller architecture is adopted for this test case.

The convolutional encoders used to compress the FOM matrices have different architectures, depending on the test case and on the value of $r$. The architectures of the CEs and CCEs are summarized in the Tables \ref{tab:conv-arch-poisson}, \ref{tab:conv-arch-adv-dom}, and \ref{tab:conv-arch-burgers}, for the three test cases, respectively.

\begin{table}[htpb!]
    \caption{\emph{Poisson test case}. Architectures of CE architectures, in terms of kernel, stride, padding, and activation function for each layer.}
    \label{tab:conv-arch-poisson}
    \centering
    \begin{tabular}{
    p{0.18\linewidth}
    p{0.13\linewidth} p{0.13\linewidth}
    p{0.07\linewidth}
    p{0.07\linewidth}
    p{0.08\linewidth}
    p{0.12\linewidth}
    }
    \toprule
    \textbf{Layer}
    &\textbf{Input dim}
    &\textbf{Output dim}
    &\textbf{Kernel size}
    &\textbf{Stride}
    &\textbf{Padding}
    &\textbf{Activation}\\
    \midrule
    Conv1 & $(1, 638, 638)$ & $(1, 127, 127)$  & $(10, 10)$ & $(5, 5)$&$(1, 1)$&Softplus\\ \midrule
    Conv2 & $(1, 127, 127)$ & $(1, 42, 42)$ & $(5, 5)$ & $(3, 3)$&$(1, 1)$&Softplus\\ \midrule
    Conv3 & $(1, 42, 42)$ & $(1, 20, 20)$ & $(5, 5)$ & $(2, 2)$&$(1, 1)$&Softplus\\ \midrule
    Conv4 ($r=2$) & $(1, 20, 20)$ & $(1, 7, 7)$ & $(4, 4)$ & $(3, 3)$&$(1, 1)$&Softplus\\
    Conv5 ($r=2$) & $(1, 7, 7)$ & $(1, 4, 4)$  & $(3, 3)$ & $(2, 2)$&$(1, 1)$&Softplus\\
    Conv6 ($r=2$) & $(1, 4, 4)$ & $(1, 2, 2)$ & $(3, 3)$ & $(2, 2)$&$(1, 1)$&\xmark \\ \midrule
    Conv4 ($r=3$) & $(1, 20, 20)$ & $(1, 10, 10)$  & $(3, 3)$ & $(2, 2)$&$(1, 1)$&Softplus\\
    Conv5 ($r=3$) & $(1, 10, 10)$ & $(1, 5, 5)$ & $(3, 3)$ & $(2, 2)$&$(1, 1)$&Softplus\\
    Conv6 ($r=3$) & $(1, 5, 5)$ & $(1, 3, 3)$  & $(3, 3)$ & $(2, 2)$&$(1, 1)$&\xmark \\ \midrule
    Conv4 ($r=4$) & $(1, 20, 20)$ & $(1, 7, 7)$ & $(4, 4)$ & $(3, 3)$&$(1, 1)$&Softplus\\
    Conv5 ($r=4$) & $(1, 7, 7)$ & $(1, 4, 4)$  & $(3, 3)$ & $(2, 2)$&$(1, 1)$&\xmark \\ \bottomrule
    \end{tabular}
\end{table}

\begin{table}[htpb!]
    \caption{\emph{Advection-diffusion test case}. Architectures of CE architectures, in terms of kernel, stride, padding, and activation function for each layer.}
    \label{tab:conv-arch-adv-dom}
    \centering
    \begin{tabular}{
    p{0.20\linewidth}
    p{0.13\linewidth} p{0.13\linewidth}
    p{0.07\linewidth}
    p{0.07\linewidth}
    p{0.08\linewidth}
    p{0.12\linewidth}
    }
    \toprule
    \textbf{Layer}
    &\textbf{Input dim}
    &\textbf{Output dim}
    &\textbf{Kernel size}
    &\textbf{Stride}
    &\textbf{Padding}
    &\textbf{Activation}\\
    \midrule
    Conv1 & $(1, 2500, 2500)$ & $(1, 623, 623)$  & $(10, 10)$ & $(4, 4)$&$(0, 0)$&Softplus\\ \midrule
    MaxPool1 & $(1, 623, 623)$ & $(1, 154, 154)$  & $(8, 8)$ & $(4, 4)$& \xmark & \xmark\\
    \midrule
    Conv2 & $(1, 154, 154)$ & $(1, 50, 50)$ & $(5, 5)$ & $(3, 3)$&$(0,0)$&Softplus\\
    \midrule
    MaxPool2 & $(1, 50, 50)$ & $(1, 16, 16)$  & $(3, 3)$ & $(3, 3)$& \xmark & \xmark\\
    \midrule
    Conv3 ($r=2$) & $(1, 16, 16)$ & $(1, 7, 7)$ & $(3, 3)$ & $(2, 2)$&$(0,0)$& \xmark \\
    MaxPool3 ($r=2$) & $(1, 7, 7)$ & $(1, 2, 2)$  & $(4, 4)$ & $(2, 2)$& \xmark & \xmark\\
    \midrule
    Conv3 ($r=3$) & $(1, 16, 16)$ & $(1, 7, 7)$ & $(3, 3)$ & $(2, 2)$&$(0,0)$& \xmark \\
    MaxPool3 ($r=3$) & $(1, 7, 7)$ & $(1, 3, 3)$  & $(3, 3)$ & $(2, 2)$& \xmark & \xmark\\
    \midrule
    Conv3 ($r=4$) & $(1, 16, 16)$ & $(1, 14, 14)$ & $(3, 3)$ & $(1, 1)$&$(0,0)$& \xmark \\
    MaxPool3 ($r=4$) & $(1, 14, 14)$ & $(1, 4, 4)$  & $(4, 4)$ & $(3, 3)$& \xmark & \xmark\\
    \bottomrule
    \end{tabular}
\end{table}

The proposed architectures include a sequence of 2D convolutional layers, activation functions, max pooling operations, and batch normalization. A continuous convolutional layer is considered in the Burgers test case (Table \ref{tab:conv-arch-burgers}) to efficiently process sparse operators. Here, $N_b$ denotes the batch size. In the previous test cases, $N_b=1$ since no parametrized operators were considered. The quantity $N_{\text{nz}}$ represents the number of non-zero entries of the operator. As described in Section~\ref{subsec:cceae}, the operator is represented using three input columns corresponding to the row indices, column indices, and values of the non-zero entries.
The filter of the continuous convolutional layer ContConv1 is a deep fully-connected neural network of input dimension $2$ (row and column indexes), hidden layers with neurons $[10, 10, 10]$, and output dimension $1$ (the continuous filter value).

\begin{table}[htpb!]
    \caption{\emph{Burgers test case}. Architectures of CCE architectures, in terms of kernel, stride, padding, and activation function for each layer.}
    \label{tab:conv-arch-burgers}
    \centering
    \begin{tabular}{
    p{0.16\linewidth}
    p{0.13\linewidth} p{0.13\linewidth}
    p{0.12\linewidth}
    p{0.07\linewidth}
    p{0.08\linewidth}
    p{0.12\linewidth}
    }
    \toprule
    \textbf{Layer}
    &\textbf{Input dim}
    &\textbf{Output dim}
    &\textbf{Kernel size}
    &\textbf{Stride}
    &\textbf{Padding}
    &\textbf{Activation}\\
    \midrule
    ContConv1 & $(N_b, 1, N_{\text{nz}}, 3)$ & $(N_b, 1, 25, 3)$  & $(4890, 4890)$ & \xmark & \xmark &\xmark\\ \midrule
    Reshape & $(N_b, 1, 25, 3)$ & $(N_b, 3, 5, 5)$ & \xmark & \xmark&\xmark&\xmark\\
    \midrule
    Conv1 ($r=2$) & $(N_b, 3, 5, 5)$  & $(N_b, 1, 2, 2)$ & $(3, 3)$ & $(2, 2)$&$(0, 0)$&Softplus\\
    BatchNorm1 ($r=2$) & $(N_b, 1, 2, 2)$  & $(N_b, 1, 2, 2)$ & \xmark & \xmark&\xmark&\xmark\\
    \midrule
    MaxPool1 ($r=3$) & $(N_b, 3, 5, 5)$  & $(N_b, 3, 4, 4)$ & $(2, 2)$ & $(1, 1)$&\xmark&\xmark\\
    Conv1 ($r=3$) & $(N_b, 3, 4, 4)$  & $(N_b, 1, 3, 3)$ & $(2, 2)$ & $(1, 1)$&$(0, 0)$&Softplus\\
    BatchNorm1 ($r=3$) & $(N_b, 1, 3, 3)$  & $(N_b, 1, 3, 3)$ & \xmark & \xmark&\xmark&\xmark
    \\ \bottomrule
    \end{tabular}
\end{table}

All the neural networks are trained with a learning rate $\num{1e-3}$, and using an early stopping based on the validation loss value. In all the test cases, a weight decay regularization with coefficient $\num{1e-2}$ was employed to improve the generalization capability of the model.

Finally, we present in Table \ref{tab:times-cae} the computational cost needed by the novel CE-AE, compared with the cost of a standard AE with the same architecture employed for the AE in the novel model.
\begin{table}[htpb!]
    \caption{Training times, number of parameters, and accuracy in terms of average test error, of the AE, CE-AE, S-CE-AE (for the second test case), and CCE-AE models. The test errors are computed by interpolating the latent variables with an RBF approach.}
    \label{tab:times-cae}
    \centering
    \begin{tabular}{
    p{0.04\linewidth}
    p{0.1\linewidth}
    p{0.06\linewidth}
    p{0.06\linewidth}
    p{0.06\linewidth}
    p{0.06\linewidth}
    p{0.06\linewidth}
    p{0.06\linewidth}
    p{0.06\linewidth}
    p{0.06\linewidth}
    p{0.06\linewidth}
    }
    \toprule
    \textbf{Test case}&\textbf{Model}&\multicolumn{3}{c}{\textbf{Average training time}}&\multicolumn{3}{c}{\textbf{Number of parameters}}&\multicolumn{3}{c}{\textbf{Average $L^2$ test error}}\\
    \midrule
    && $r=2$&$r=3$&$r=4$& $r=2$&$r=3$&$r=4$& $r=2$&$r=3$&$r=4$\\
    \cmidrule{3-5} \cmidrule{6-8} \cmidrule{9-11}
    \multirow{3}{*}{1}&AE&$24.6$ &$381$ & $370$& \multicolumn{3}{c}{$\sim134$K} & $0.152$ & {$\mathbf{0.026}$} & $0.026$\\
    &CE-AE&$221$ &$606$ &$116$ & {$+190$} &{$+183$}& {$+180$} & $\mathbf{0.041}$ & $0.035$ & {$\mathbf{0.025}$} \\
    &\small S-CE-AE& $245$ & $527$& $109$& {$+190$} &{$+183$} &{$+180$} & $0.045$ & $0.032$ & $0.028$ \\
    \midrule
    && $r=2$&$r=3$&$r=4$& $r=2$&$r=3$&$r=4$& $r=2$&$r=3$&$r=4$\\ \cmidrule{3-5} \cmidrule{6-8} \cmidrule{9-11}
    \multirow{2}{*}{2}&AE& $920$ & $999$ & $1078$ &\multicolumn{3}{c}{$\sim516$K} & $0.032$ & $0.024$ & $0.010$ \\
    &CE-AE& $1331$ & $1461$ & $1882$ & $+219$ & $+212$ & $+219$ & $\mathbf{0.0035}$ & $\mathbf{0.0082}$ & $\mathbf{0.0030}$ \\
    \midrule
    && $r=2$&$r=3$&& $r=2$&$r=3$&& $r=2$&$r=3$&\\ \cmidrule{3-4} \cmidrule{6-7} \cmidrule{9-10}
    \multirow{2}{*}{3}&AE& $55.1$ & $39$ & &\multicolumn{2}{c}{$\sim513$K} & & $0.50$ & $0.52$ & \\
    &CCE-AE&$9889$ & $3454$ & & $+399$ & $+384$ &  & $\mathbf{0.085}$ & $\mathbf{0.084}$ & \\
    \bottomrule
    \end{tabular}
\end{table}
Regarding the computational time, for the first two test cases the proposed models exhibit a slightly higher computational cost compared to the AE, while remaining of the same order of magnitude. This increase is mainly due to the additional convolutional encoder, which introduces only a limited number of trainable parameters compared to the AE. Therefore, the increase in model complexity remains marginal.

Conversely, for the Burgers' case, the computational cost is higher due to two main factors: the presence of the continuous convolutional layer, which involves more expensive routines (e.g., the filter mapping function), and the time-loop integration within the training procedure, which requires the reduced dynamics to be evaluated at each epoch. Nevertheless, also in this case, the number of additional trainable parameters introduced with respect to the AE remains negligible.

Despite the additional computational effort, the proposed approaches achieve improved accuracy compared to the AE in the testing phase, except for the first test case with $r=3$.

\section{Conclusions and Future Perspectives}
\label{sec:conclusions}

In this work we introduced a novel perspective on \textbf{nonlinear model order reduction} based on the compression of full-order operators through convolutional neural network architectures. Unlike most existing nonlinear ROM methodologies, which focus on learning nonlinear representations of the solution manifold, the proposed framework aims at learning compact representations of the operators defining the governing equations. This preserves the algebraic structure of the reduced problem, allowing the reduced system to be solved through standard numerical techniques while extending classical projection-based ROMs beyond linear or quadratic subspaces.

The proposed methodology was assessed on three benchmark problems of increasing complexity. The numerical experiments demonstrated that the operator-compression strategy is capable of producing accurate ROMs while maintaining a compact latent representation. Moreover, the results highlighted that incorporating the solution of the reduced system directly into the training process improves the quality of the learned latent space with respect to a standard fully-connected AE.

Interestingly, even without explicitly enforcing structural constraints, the proposed architectures naturally learned reduced operators that reflected the algebraic properties of the corresponding full-order operators. In the Poisson test case, symmetric full-order operators consistently led to nearly symmetric reduced operators, indicating that these structural features can emerge naturally during training.

Furthermore, incorporating the reduced-system solution into the training objective acts as an effective \textbf{stabilization mechanism} for the autoencoder optimization. Although it introduces an additional loss contribution, it guides the optimization towards more informative latent representations, leading to improved training stability and, in several cases, to lower reconstruction losses than those achieved by a standard autoencoder.

Finally, the introduction of continuous convolutional layers provided an effective strategy for exploiting the sparsity of differential operators, making the approach scalable to large discretized systems and naturally applicable to operators defined on unstructured meshes.

Although the proposed framework should be regarded as a preliminary step towards nonlinear operator-based ROMs, the obtained results suggest several promising research directions. A first natural extension consists in applying the methodology to \textbf{more realistic and large-scale CFD problems}, including the incompressible Navier--Stokes equations and multiphysics applications. Furthermore, more expressive parameter-to-latent mappings, such as neural networks, could replace the interpolation strategies adopted in this work, improving the generalization capability over complex parameter spaces.

Finally,\textbf{ combining continuous convolutional networks for both the operators and the solution fields} would provide a fully discretization-independent framework, enabling nonlinear ROMs that naturally generalize across different computational grids.

\bibliographystyle{abbrv}
\bibliography{main}
\end{document}